\documentclass[a4paper,10pt,oneside,openany]{article}
\usepackage{graphics}
\usepackage{float}
\usepackage[dvips]{graphicx}
\usepackage{subfigure}
\usepackage{latexsym,euscript,epic,eepic,psfig}
\newtheorem{theo}{Theorem}
\newtheorem{rem}{Remark}

\newtheorem{cor}{Corollary}
\newtheorem{lem}{Lemma}
\newcommand{\cvP}{\stackrel{{\rm{P}}}{\longrightarrow}}
\newcommand{\cvps}{\stackrel{{\rm{a.s.}}}{\longrightarrow}}
\newcommand{\cvloi}{\stackrel{d}{\longrightarrow}}
\newcommand{\egloi}{\stackrel{d}{=}}
\newcommand{\eqproba}{\stackrel{{\rm{P}}}{\sim}}
\newcommand{\R}{{\rm{I\!\!R}}}
\newcommand{\N}{{\rm{I\!\!N}}}
\newcommand{\ind}{{\rm{1\!\!I}}}
\newcommand{\proba}{{\rm{P}}}
\newcommand{\findemo}{\hspace*{\fill} $\spadesuit$}

\begin{document}

\begin{center}
{\large{\bf{
A Pickands type estimator of the extreme value index \\
}}}
Laurent Gardes and St\'{e}phane Girard {\footnote{{\tt{Corresponding author: Stephane.Girard@imag.fr}}}}
\end{center}
\vspace*{2mm}\par
\begin{center}
SMS/LMC, Universit\'{e} de Grenoble 1, BP 53, 38 041 Grenoble Cedex 9, France.
\end{center}
\vspace*{2mm}\par
{\bf{Abstract}} $-$ One of the main goals of extreme value analysis is
to estimate the probability of rare events given a sample from an
unknown distribution. The upper tail behavior of this distribution is
described by the extreme value index. We present a new estimator of
the extreme value index adapted to any domain of attraction. Its
construction is similar to the one of Pickands' estimator. Its weak
consistency and its asymptotic distribution are established and a bias reduction method is proposed. Our estimator is compared with classical extreme value index estimators through a simulation study. \\
\vspace*{8mm}\par

\section{Introduction}

Suppose one is given a sequence $X_1,\ldots,X_n$ of independent and
identically distributed (i.i.d.) observations from some distribution
function $F$. Suppose there exist sequences $a_n > 0$ and $b_n$ and
some $\xi \in \R$ such that:
\begin{equation}
\label{loimax}
\lim_{n \to \infty} \proba \left [ \frac{ \max(X_1,\ldots,X_n) -
b_n}{a_n} \leq x \right ] = G_{\xi}(x), 
\end{equation}
with $G_{\xi}(x)=\exp[-(1+\xi x)_{+}^{-1/\xi}]$ if $\xi \neq 0$ and
$G_0(x)=\exp[-{\rm{e}}^{-x}]$, where $y_+=\max(0,y)$. Necessary and
sufficient conditions on $F$ for the convergence (\ref{loimax}) to the
extreme value distribution $G_{\xi}$ can be found in \cite{gne43}. The aim of this
paper is the definition of a new estimator of the extreme value index $\xi \in
\R$. This parameter drives the decay of
the tail distribution: as a power function if $\xi > 0$ (Pareto,
Burr, Student's, Log-gamma distributions, etc \ldots), exponentially if $\xi=0$
(Exponential, Normal, Log-normal, Gamma distributions, etc \ldots) and
with finite right endpoint if $\xi <0$ (Uniform, Beta, Reversed
Pareto, Reversed Burr distributions, etc \ldots). The knowledge of
$\xi$ is for example of high interest for extreme quantile estimation
which arises in a lot of applications \cite{embklumik97} such as
finance, insurance, hydrology, etc  \ldots \ There is a substantial number of publications dedicated to the estimation of this extreme value index, especially on the heavy tailed
distribution context ($\xi > 0$) (see Beirlant et {\it{al.}} \cite{beivynteu96}, Feueverger and Hall \cite{feuhal99} and, for a recent overview of this literature, see Cs\"{o}rgo and Viharos \cite{csovih98}). The most popular estimator in this case is the Hill estimator \cite{hil75} defined by:
\[ {\hat{\xi}}_{k,n}^{H} = \frac{1}{k} \sum_{i=1}^{k} \ln(X_{n-i+1,n}) -
\ln(X_{n-k,n}), \ {\rm{for}} \ k=1,\ldots,n-1, \]
where $X_{1,n} \leq \ldots \leq X_{n,n}$ correspond to the random variables $X_1,\ldots,X_n$ rearranged in ascending order. The consistency and
the asymptotic normality of this estimator are proved for example by
Davis and Resnick \cite{davres84}, Cs\"{o}rgo and Mason
\cite{csomas85}, etc \ldots \\
The general case $\xi \in \R$ has been less extensively studied.
Dekkers, Einmahl and de Haan \cite{dekeinhaa89} have adapted the
estimator proposed by Hill to this situation. Another estimator was proposed by Pickands \cite{pic75}:
\[ {\hat{\xi}}_{k,n}^{P} = \frac{1}{\ln(2)} \ln \left (
\frac{X_{n-k+1,n}-X_{n-2k+1,n}}{X_{n-2k+1,n}-X_{n-4k+1,n}} \right ), \
{\rm{for}} \ k=1,\ldots, \lfloor n/4 \rfloor, \]
where $\lfloor x \rfloor$ denotes the integer part of $x$. Weak and strong
consistency as well as asymptotic normality of ${\hat{\xi}}_{k,n}^{P}$
were established by Dekkers and de Haan \cite{dekhaa89}. Proofs are
based on the following well known result: Let $U$ be the tail
quantile function of the distribution function $F$ defined by
\[ U(x) = \left( \frac{1}{1-F(x)} \right )^{\leftarrow} \]
(the arrow means inverse function) and let
\[ \varphi_{t}(x)=\int_1^xu^{t-1}du, \ x>0, \ t \in \R. \]
Relation (\ref{loimax}) holds if
and only if (see de Haan \cite{haa84}), there exist a positive measurable function $a$ such that uniformly locally on $x >0$,
\begin{equation}
\label{equiv1loimax}
\lim_{t \to \infty} \frac{U(tx)-U(t)}{a(t)} = \varphi_{\xi}(x).
\end{equation}
Clearly, relation (\ref{equiv1loimax}) implies
that uniformly locally on $x,y>0$, $y \neq 1$,
\begin{equation}
\label{equiv2loimax}
\lim_{t \to \infty} \frac{U(tx)-U(t)}{U(ty)-U(t)} =
\frac{\varphi_{\xi}(x)}{\varphi_{\xi}(y)}.
\end{equation}
Thus, by substituting in (\ref{equiv2loimax}) $U$ by
${\hat{U}}_n=[1/(1-{\hat{F}}_n)]^{\leftarrow}$ (${\hat{F}}_n$ denoting
the empirical distribution function), $t$ by $n/(2k)$, $x$ by $1/2$
and $y$ by $2$, and remarking that ${\hat{U}}_n(n/k)=X_{n-k+1,n}$, we
have asymptotically
\begin{equation}
\label{picequation}
2^{-\xi} \frac{X_{n-k+1,n}-X_{n-2k+1,n}}{X_{n-2k+1,n}-X_{n-4k+1,n}}
= 1.
\end{equation}
Pickands' estimator ${\hat{\xi}}_{k,n}^{P}$ is the solution of the equation
(\ref{picequation}). One can notice that this estimator does not take into account of the extreme observations $X_{n-k+2,n}, \ldots,X_{n,n}$.\\
 \ \\
In the next section, we define a new estimator of the extreme value
index $\xi$ when $\xi \in \R$. This estimator is similar to the one of
Pickands but exploiting the information given by the spacing between
$X_{n-k+1,n}$ and $X_{n,n}$. Weak consistency and asymptotic
distribution are established in section~3 and a bias corrected
estimator is introduced. Section~4 is devoted to the proofs of the main
results and a simulation study is presented in section~5.

\section{Estimation of the extreme value index}

We propose to estimate the extreme value index $\xi \in \R$ by
${\hat{\xi}}_{k,n}$ defined as the root of the equation
in $\theta$:
\begin{equation}
\label{estimatordef}
\left \{  \frac{\varphi_{\theta}(1/k')}{\varphi_{\theta}(1/k)} \right \}
\frac{X_{n-k+1,n}-X_{n,n}}{X_{n-k'+1,n}-X_{n,n}} = 1, \ {\rm{for}} \ 1
< k' < k < n.
\end{equation}
We can show (see Gardes \cite{gar03}, Appendix B) that (\ref{estimatordef}) admits
an unique solution. This estimator applies to all real $\xi$ and, as
Pickands' estimator, remains unaffected when the scale or location of
the data are changed. Furthermore, as we will see on a simulation
study, the behavior of ${\hat{\xi}}_{k,n}$ is less influenced by the
parameter $k$ than Pickands' estimator. One can justify the
definition of ${\hat{\xi}}_{k,n}$ by the two following lemmas: 

\begin{lem}
\label{lem1}
{\rm{ Suppose that relation (\ref{loimax}) holds. Then, if $\xi < 0$,
\[ \lim_{t \to \infty, \ x \to \infty} \frac{U(tx)-U(t)}{a(t)} = -\frac{1}{\xi} \]
and if $\xi \geq 0$,
\[ \lim_{t \to \infty, \ x \to \infty} \frac{U(tx)-U(t)}{a(t)} = +\infty. \]
}}
\end{lem}
 \ \\
As a consequence of Lemma \ref{lem1}, we have: 

\begin{lem}
\label{lem2}
{\rm{Suppose that relation (\ref{loimax}) holds. Then,
\begin{equation}
\label{justiesti}
\left \{ \frac{\varphi_{\xi}(y)}{\varphi_{\xi}(x)} \right \}
\frac{U(tx)-U(t)}{U(ty)-U(t)} \to  1,
\end{equation}
as $t \to \infty$, $x \to 0$ with $ty \to \infty$ and $x/y \to d>0$.
}}
\end{lem}
 \ \\
Lemma \ref{lem1} and Lemma \ref{lem2} can be seen as an extension of
respectively (\ref{equiv1loimax}) and (\ref{equiv2loimax}) when $x$
and $y$ are going to zero or infinity. The proofs of Lemma~\ref{lem1} and Lemma~\ref{lem2} are postponed to the Appendix.
By substituting in (\ref{justiesti}) $U$ by ${\hat{U}}_n$, $t$ by $n$,
$x$ by $1/k'(n)=1/k'$ and $y$ by $1/k(n)=1/k$ with $k/k' \to
c > 1$, $k \to \infty$ and $k/n \to 0$ as $n \to \infty$, we have asymptotically:
\[ \left \{ \frac{\varphi_{\xi}(1/k')}{\varphi_{\xi}(1/k)} \right \}
\frac{X_{n-k+1,n}-X_{n,n}}{X_{n-k'+1,n}-X_{n,n}} = 1, \]
which is an intuitive justification for the definition of
${\hat{\xi}}_{k,n}$. The next section is dedicated to the study of
${\hat{\xi}}_{k,n}$ asymptotical properties.

\section{Main results}

\subsection{Asymptotic properties}

We first state the weak consistency of ${\hat{\xi}}_{k,n}$ under some conditions on $k$ and $k'$.

\begin{theo}
\label{cvproba}
{\rm{Suppose that relation (\ref{loimax}) holds. If $k/k' \to c>1$, $k
\to \infty$ and $k/n \to 0$ as $n \to \infty$, then ${\hat{\xi}}_{k,n}
\cvP \xi$.}}
\end{theo}

\begin{rem}
{\rm{Similar conditions on $k$ are used by Dekkers and de Haan \cite{dekhaa89} to prove that Pickands' estimator ${\hat{\xi}}_{k,n}^{P}$ is weakly consistent.}}
\end{rem}

\noindent To establish the asymptotic distribution of the estimator
${\hat{\xi}}_{k,n}$, additionnary conditions are introduced. The first
of them is a cornerstone in all proofs of asymptotic normality for
extreme value estimators. \\
\begin{itemize}
\item[$(H1)$ $-$] $U$ has a positive derivative and there exist a slowly
varying function $\ell$ such that $U'(x)=x^{\xi-1}\ell(x)$. 
\end{itemize}
We refer to \cite{bingolteu87} for more details on slow variation
theory. The next condition controls the uniform rate of convergence of $\ell(tx)/\ell(x)$ to $1$ as $x \to \infty$. Let $\delta=\min(-\xi,1/2)$ and introduce the random variables
$K_{k,n}={\bar{F}}(X_{n-k+1,n})/{\bar{F}}(X_{n,n})$ and
$N_n=1/{\bar{F}}(X_{n,n})$ where ${\bar{F}}$ is the survival function (${\bar{F}}=1-F$). \\
\begin{itemize}
\item[$(H2)$ $-$] 
\[ \varphi_{\delta}(k') \sup_{t \in [1,K_{k',n}]} \left | \frac{\ell(tN_n/K_{k',n})}{\ell(N_n/K_{k',n})} - 1 \right | \cvP 0. \]
\end{itemize}
 \ \\
Our second main result is the following:

\begin{theo}
\label{cvloi}
{\rm{Let $V_k(\xi)=\varphi_{\delta}(k)[(\ln(k)-1)\ind\{\xi \geq
0\}+1]$. Under the conditions of Theorem \ref{cvproba} (with $k=ck'$) and if $(H1)$ and $(H2)$ are satisfied, we have for all $t \in \R$:
\begin{equation}
\label{loiasymp}
\lim_{n \to \infty} \proba [V_k(\xi)({\hat{\xi}}_{k,n}-\xi) \leq t] = \left \{ \begin{array}{l l} 
\exp(-{\rm{e}}^{-t}) & {\rm{if}} \ 0<\xi, \\
\exp(-{\rm{e}}^{-t/2}) & {\rm{if}} \ \xi = 0, \\
\exp \left [ - \left [ 1+ t \ln(c) / \varphi_{\xi}(1/c) \right
]^{-1/\xi} \right ] & {\rm{if}} \ -1/2 < \xi <0, \\
\Phi \left [ - t c^{-\xi} \ln(c)/(2 \xi \sigma) \right ] & {\rm{if}} \ \xi<-1/2, \end{array} \right . 
\end{equation}
where $\sigma=c^{-\xi} (c-1)^{1/2}$ and $\Phi$ is the cumulative distribution function of the standard normal distribution.}}
\end{theo}

\begin{rem}
{\rm{ \ \
\begin{itemize} 
\item[i)] Theorem \ref{cvloi} states that the asymptotic distribution
of ${\hat{\xi}}_{k,n}$ is Gaussian if $\xi<-1/2$ and an extreme value
distribution if $\xi>-1/2$. If $\xi=-1/2$, we prove that
$V_k(\xi)({\hat{\xi}}_{k,n}-\xi)$ converges to a non-degenerate
distribution with non explicit cumulative distribution function. In
fact, as it will appear in the next section, the limit distribution of ${\hat{\xi}}_{k,n}$ is driven by $k'/K_{k',n}$ if $\xi<-1/2$, by $Y_n$ if $\xi>-1/2$ with 
\[ Y_n=\sqrt{k}/\sigma [\varphi_{\xi}(K_{k',n}/K_{k,n})-\varphi_{\xi}(1/c)] \]
and by both of them if $\xi=-1/2$.
\item[ii)] $(H1)$ and $(H2)$ are second order conditions on the tail quantile function $U$. Similar conditions are used by Dekkers and de Haan \cite{dekhaa89} to establish the asymptotic distribution of Pickands' estimator. 
\item[iii)] Let $\mu(\xi)=\gamma \ind\{\xi>0\}-[1-\Gamma(1-\xi)] \varphi_{\xi}(1/c)/\ln(c) \ind \{-1/2<\xi<0\}$, where $\gamma$ is the Euler constant and $\Gamma$ is the gamma function. Theorem \ref{cvloi} entails that $V_k(\xi)({\hat{\xi}}_{k,n}-\xi)$ converges to a distribution of mean $\mu(\xi)$ if $\xi \neq 0$ and $\xi \neq -1/2$. This suggests to define the bias corrected estimator:
\[ {\hat{\xi}}_{k,n}^{*} = {\hat{\xi}}_{k,n} - \frac{\mu({\hat{\xi}}_{k,n})}{V_k({\hat{\xi}}_{k,n})}. \]
As we will see on a simulation study (see section \ref{subsecsimul1}),
this bias correction improves the behavior of our estimator in most finite sample situations.
\end{itemize} }}
\end{rem}

\subsection{Examples}
\label{secex}

Let $\alpha$, $\beta >0$, $\theta \in \R \backslash
\{0\}$ and define $\ln_2(x)=\ln(\ln(x))$, $x>1$. The two following
models of slowly varying functions $\ell$ are considered: 
\[ \ell(x)=\alpha+\theta x^{-\beta}+o(x^{-\beta}), \ \ {\rm{\bf{(Model \
A)}}}, \]
\[ \ell(x)=\theta[\ln(x)]^{-\beta} \left \{ 1 + O \left (
\frac{\ln_2(x)}{\ln(x)} \right ) \right \}, \ \ {\rm{\bf{(Model \
B).}}} \]
Model A has been first introduced by Hall \cite{hal82}. In
both models, the parameter $\beta$ tunes the decay of the
slowly varying function $\ell$. The conditions that should be
satisfied by models A and B to insure convergence (\ref{loiasymp}) are
given in Corollary \ref{corex}. In both cases, the best rate of
convergence of ${\hat{\xi}}_{k,n}$ is also established. Some examples
of distributions satisfying the assumptions of Corollary 
\ref{corex} are presented in Table \ref{tab1}. In the sequel, the
following notation is adopted. Let $(u_n)$ and $(v_n)$ be two non negative
deterministic sequences. The notation $u_n \asymp v_n$ means that
\[ 0 < \liminf \frac{u_n}{v_n} \leq \limsup \frac{u_n}{v_n} < \infty. \]

\begin{cor}
\label{corex}
{\rm{ Suppose that $k=ck'$, $k \to \infty$ and that $F$ satisfies assumption $(H1)$ with a slowly varying function asymptotically monotone.
\begin{itemize}
\item[i)] If $\ell$ belongs to Model A and if $\varphi_{\delta}(k') \left ( n/k' \right )^{-\beta} \to 0$ then convergence (\ref{loiasymp}) holds. In this case, the best rate of convergence of ${\hat{\xi}}_{k,n}$ is given by:
\[
V_k(\xi) \asymp \left \{ 
\begin{array}{l l} 
\ln (n) & {\rm{if}} \ 0<\xi, \\
\ln^2(n) & {\rm{if}} \ \xi=0, \\
n^{\delta/[1+\delta/\beta]-\varepsilon} & {\rm{if}} \ \xi < 0, \\
\end{array} \right . 
\]
where $\varepsilon \in ]0,\delta/(1+\delta/\beta)[$ is arbitrarly small.
\item[ii)] If $\ell$ belongs to Model B and if $\varphi_{\delta}(k') \ln(k')/\ln(n) \to 0$ and $\varphi_{\delta}(k')
\ln_2(n)/\ln(n) \to 0$ then, convergence (\ref{loiasymp}) holds. Furthermore, the best rate of convergence of ${\hat{\xi}}_{k,n}$ is given by:
\[
V_k(\xi) \asymp \left \{ 
\begin{array}{l l} 
\ln_2 (n) & {\rm{if}} \ 0<\xi, \\ 
\ln_2^2(n) & {\rm{if}} \ \xi=0, \\
\ln^{1-\varepsilon}(n) & {\rm{if}} \ \xi < 0, \\ 
\end{array} \right .
\]
where $\varepsilon \in ]0,1[$ is arbitrarly small.
\end{itemize}
}} \end{cor}

\begin{rem}
{\rm{
This corollary points out the fact that the case $\xi<0$ is more
favorable to our estimator, i.e. its convergence is faster than in the
case $\xi \geq 0$. This is illustrated by the simulation study (see section~\ref{secsimul}). \\

}}
\end{rem}

\begin{table}
\begin{center}
\begin{tabular}{|c|c|c|c|c|}
\hline 
Distribution & Cumulative distribution  & Model & $\beta$ &
Best rate of  \\
 \ & function & \ & \ & convergence  \\
\hline 
\hline
Weibull$_{\rm{M}}$  & $\exp[-(1+\xi x)^{-1/\xi}]$, & A & 1 & $\ln(n)$ if $0<\xi$ \\
 ($\xi \in \R$) & for $x$ such that $1+\xi x >0$. & \ & \ & $\ln^2(n)$ if $\xi=0$ \\
 \ & \ & \ & \ & $n^{\delta/(1+\delta/\beta)-\varepsilon}$ if $\xi<0$ \\
\hline
Burr  & $1-[w/(w+x^{\tau})]^{\lambda}$, for $x>0$, & A & $1/\lambda$ & $\ln(n)$ \\
 ($\xi>0$) & with $w$, $\lambda$, $\tau >0$, $\xi=1/(\lambda \tau)$, & \ & \ & \ \\
\hline
Fr\'{e}chet  & $\exp(-x^{-1/\xi})$ for $x>0$, & A & 1 &
$\ln(n)$ \\
($\xi >0$) & \ & \ & \ & \ \\
\hline 
Weibull & $1-\exp(-\lambda x^{\tau})$, for $x>0$, & B & $1-1/\tau$ & $\ln_2^2(n)$ \\
 ($\xi=0$) & with $\lambda$, $\tau > 0$. & \ & \ & \ \\
\hline
Normal & $\int_{-\infty}^x 1/\sqrt{2\pi} {\rm{e}}^{-t^2/2} dt$ & B &
1/2 & $\ln_2^2(n)$. \\
($\xi=0$) & \ & \ & \ & \ \\
\hline
Reversed Burr & $1-[w/(w+(x_F-x)^{-\tau})]^{\lambda}$, & A & $1/\lambda$ & \ \\
($\xi<0$) \ & for $x < x_F$, with $\xi=1/(\lambda \tau)$. & \ & \ & $n^{\delta/(1+\delta/\beta)-\varepsilon}$ \\
 \ & $w$, $\lambda$, $\tau >0$, $x_F \in \R$. & \ & \ & \ \\
\hline
\end{tabular}
\end{center}
\caption{Examples of distributions satisfying the assumptions of
Corollary \ref{corex}}
\label{tab1}
\end{table}

\section{Proofs of the main results}

This section is devoted to the proof of Theorem \ref{cvproba} and Theorem
\ref{cvloi}. Proofs of lemmas are postponed to the appendix.

\subsection{Preliminary results}

The following function will play an important role. Let 
\[ H_n(x)= \left \{ \frac{\varphi_{x}(1/k')}{\varphi_{x}(1/k)} \right \}
\frac{X_{n-k+1,n}-X_{n,n}}{X_{n-k'+1,n}-X_{n,n}} = \left \{ \frac{\varphi_{x}(1/k')}{\varphi_{x}(1/k)} \right \} (1+Z_n), \]
with $Z_n=(X_{n-k+1,n}-X_{n-k'+1,n})/(X_{n-k'+1,n}-X_{n,n})$.

\begin{lem}
\label{lemZ}
{\rm{Under the conditions of Theorem \ref{cvproba}, 
\begin{itemize} 
\item[i)] $Z_n \egloi [U(N_n/K_{k,n})-U(N_n/K_{k',n})]/[U(N_n/K_{k',n})-U(N_n)]$. 
\item[ii)] $N_n \cvps +\infty$, $K_{k,n} \cvps +\infty$, $K_{k,n}/K_{k',n} \cvps c$, $K_{k,n}/N_n \cvps 0$ and $Z_n \cvps \max(0,c^{-\xi}-1)$ as $n \to \infty$.
\item[iii)] $k'/K_{k',n} \cvloi {\rm{Exp}}(1)$. If moreover $k=ck'$,
then $Y_n \cvloi {\cal{N}}(0,1)$.
\end{itemize}
}}
\end{lem}

\subsection{Proof of Theorem \ref{cvproba}}

We shall need the following result: 

\begin{lem}
\label{lem3}
{\rm{Suppose that relation (\ref{loimax}) holds. If
$\xi > 0$, for all $\eta_1 \in ]0,\xi[$, $\eta_2 > 0$, there exist
$t_0$, $\beta_1,\beta_2,{\tilde{\beta_1}},{\tilde{\beta_2}}>0$ such that, for all $t \geq t_0$ and $x>0$,
\[ \beta_1 x^{\xi-\eta_1}-{\tilde{\beta_1}} \leq \frac{U(tx)-U(t)}{a(t)} \leq
\beta_2 x^{\xi+\eta_2}-{\tilde{\beta_2}}. \]
If $\xi=0$, for all $\eta > 0$, there exist $t_0$, $\beta>0$ such that,
for all $t \geq t_0$ and $x>0$,
\[ \frac{U(tx)-U(t)}{a(t)} \leq \beta x^{\eta}. \]
}}
\end{lem}

\noindent {\bf{Proof of Theorem \ref{cvproba} $-$ }} We have to show
that, for all $\varepsilon > 0$,
\begin{equation}
\label{amontrer}
\lim_{n \to \infty} \proba [ | {\hat{\xi}}_{k,n}-\xi| > \varepsilon ] = 0.
\end{equation}
Remark that if $\xi \neq 0$, proving (\ref{amontrer}) for all
$\varepsilon >0$ reduces to demonstrate (\ref{amontrer}) for all
$0<\varepsilon <|\xi|$. 
Since $H_n$ is a non-decreasing function (see \cite{gar03}, Appendix B)
and since $H_n({\hat{\xi}}_{k,n})=1$, we have:
\[ \proba[|{\hat{\xi}}_{k,n}-\xi|>\varepsilon] =
\proba[H_n(\xi+\varepsilon)<1]+\proba[H_n(\xi-\varepsilon)>1]. \]
To prove Theorem \ref{cvproba} it is sufficient to establish
that $\proba[H_n(\xi+\varepsilon)\geq 1] \to 1$ and
$\proba[H_n(\xi-\varepsilon)\leq 1] \to 1$ as $n \to \infty$.
The two following expansions hold:  
\begin{equation}
\label{deterpart}
\frac{\varphi_{t}(1/k')}{\varphi_{t}(1/k)} = 1+k^{-t}-k'^{-t}+o(k^{-t}) \ {\rm{if}}
\ t>0 \ {\rm{and}} \ \frac{\varphi_{t}(1/k')}{\varphi_{t}(1/k)} \to c^t \ {\rm{as}} \ n
\to \infty \ {\rm{if}} \ t<0. 
\end{equation}
The two following cases are considered separately:\\
{\underline{If $\xi \geq 0$}}, (\ref{deterpart}) and Lemma \ref{lemZ} ii) imply that $H_n(\xi+\varepsilon) = 1+Z_n+k^{-(\xi+\varepsilon)}-k'^{-(\xi+\varepsilon)}+o_{{\rm{P}}}[k^{-(\xi+\varepsilon)}]$ since $\xi+\varepsilon>0$. Furthermore, from Lemma \ref{lemZ} i),
\begin{equation}
\label{randompart}
Z_n \egloi {\underbrace{\frac{U(N_n/K_{k,n})-U(N_n/K_{k',n})}{a(N_n/K_{k',n})}}_{Z_{1,n}}} \times {\underbrace{\frac{a(N_n/K_{k',n})}{U(N_n/K_{k',n})-U(N_n)}}_{Z_{2,n}}}, 
\end{equation}
with 
\begin{equation}
\label{aa}
Z_{1,n} \cvps \varphi_{\xi}(c^{-1}),
\end{equation}
from (\ref{equiv1loimax}) and Lemma \ref{lemZ} ii) and
\begin{equation}
\label{bb}
\beta_1 K_{k',n}^{\xi-\eta_1}-{\tilde{\beta_1}} \leq -\frac{1}{Z_{2,n}} \leq \beta_2 K_{k',n}^{\xi+\eta_2}-{\tilde{\beta_2}},
\end{equation}
from Lemma \ref{lem3}. Since, from Lemma \ref{lemZ} iii), $k'/K_{k',n}$ converges to a standard exponential distribution, we deduce from (\ref{randompart})-(\ref{bb}) that $k'^{-(\xi+\varepsilon)}/Z_n \cvP 0$ which entails that $H_n(\xi+\varepsilon)-1 \eqproba Z_n>0$ i.e. that $\proba[H_n(\xi+\varepsilon) \geq 1] \to 1$. Similarly, we prove that $\proba[H_n(\xi-\varepsilon) \leq 1] \to 1$.\\
{\underline{If $\xi<0$}}, expansions (\ref{deterpart}) and Lemma \ref{lemZ} ii)
imply that $H_n(\xi+\varepsilon) \to c^{\varepsilon}>1$ since $\xi +
\varepsilon >0$. Thus, $\proba[H_n(\xi+\varepsilon) \geq 1] \to 1$. In
the same way, we prove that $\proba[H_n(\xi-\varepsilon) \leq 1] \to 1$. \findemo

\subsection{Proof of Theorem \ref{cvloi}}

Let us define the function:
\[ \varphi_t^{*}(x)= \left \{ \begin{array}{l l}
1+tx & {\rm{if}} \ t \neq 0, \\
{\rm{e}}^x & {\rm{if}} \ t=0. \end{array} \right . \]
To prove Theorem \ref{cvloi}, two auxiliary results are necessary.
Lemma \ref{lemphi} is dedicated to the study of the function $\varphi_{t}^*$.

\begin{lem}
\label{lemphi}
{\rm{ \ \ 
\begin{itemize}
\item[I)] For $x \in (0,\infty)$,
$\varphi_t^*[\varphi_t(x)]=x^{t+\ind\{t=0\}}$.
\end{itemize}
Let $(u_n)$ and $(v_n)$ be two sequences such that $u_n \sim v_n$
(i.e. $u_n/v_n \to 1$).
\begin{itemize}
\item[II)] Let $t=0$. If $u_n \to \infty$ and $u_n-v_n \to \alpha$ then $\varphi_t^*(v_n) \sim \varphi_t^{*}(u_n){\rm{e}}^{-\alpha}$.
\item[III)] Let $t \neq 0$. If $u_n \to \infty$ then $\varphi_t^*(v_n) \sim \varphi_t^{*}(u_n)$.
\item[IV)] Let $t \neq 0$. If $u_n \to -1/t$ with
$v_n=u_n(1+\varepsilon_n)$, then: 
\begin{itemize}
\item[i)] If moreover $\varepsilon_n/\varphi_t^*(u_n) \sim \alpha_n$
where $\alpha_n$ does not converge to $\infty$ or to $1$, then \\ 
$\varphi_t^*(v_n) \sim \varphi_t^{*}(u_n)(1-\alpha_n)$.
\item[ii)] If moreover $\varepsilon_n/\varphi_t^*(u_n) \to \infty$ then $\varphi_t^*(v_n) \sim -\varepsilon_n$.
\end{itemize}
\end{itemize}
}}
\end{lem}

\noindent The proof of this basic result is not detailed here. Clearly,
the distribution of $H_n(x)$ is determined by $Z_n$. The following
lemma provides the asymptotic distribution of $Z_n$.

\begin{lem}
\label{lemcvloiZ}
{\rm{ Under the conditions of Theorem \ref{cvloi},
\[ \lim_{n \to \infty} \proba \left [ k'^{\delta-\ind\{\xi=0\}} \varphi_{\xi}^* \left ( - \frac{\varphi_{\xi}(1/c)}{Z_n} \right ) \leq t \right ] = \left \{ \begin{array}{l l}
\exp(-t^{-1/\xi}) & {\rm{if}} \ 0<\xi,\\
\exp(-t^{-1}) & {\rm{if}} \ \xi =0, \\
1-\exp(-t^{-1/\xi}) & {\rm{if}} \ -1/2<\xi<0, \\
\proba[T<t\sqrt{c}] & {\rm{if}} \ \xi=-1/2, \\
\Phi \left [ -t \varphi_{\xi}(1/c)\sqrt{c}/\sigma \right ]  &
{\rm{if}} \ \xi<-1/2, \end{array} \right . \]
where the random variable $T$ is defined as the limit in distribution
of
\[ T_n =
\sqrt{\frac{k}{K_{k',n}}}+\frac{\sigma}{\varphi_{\xi}(1/c)}Y_n, \]
which is non-degenerate from Lemma \ref{lemZ} iii).
}}
\end{lem}
 \ \\
\noindent {\bf{Proof of Theorem \ref{cvloi} $-$ }} Let $F_n(t)=\proba[V_k(\xi)({\hat{\xi}}_{k,n}-\xi) \leq t]$. We have,
\[ F_n(t)=\proba \left [ {\hat{\xi}}_{k,n} \leq \xi + t/V_k(\xi) \right ] = \proba [H_n(\xi+t/V_k(\xi)) \geq 1]=\proba \left [ (1+Z_n) \frac{\varphi_{-\xi-t/V_k(\xi)}(k/c)}{\varphi_{-\xi-t/V_k(\xi)}(k)} \geq 1 \right ], \]
since $H_n$ is a non-decreasing function and since
$H_n({\hat{\xi}}_{k,n})=1$.
%Remarking that, for $n$ large enough and for all $\xi \in \R$,
%\[ \frac{\varphi_{-\xi-t/V_k(\xi)}(k/c)}{\varphi_{-\xi-t/V_k(\xi)}(k)} > 0, \]
%implies that
%\[ F_n(t)=\proba \left [ Z_n \geq \frac{\varphi_{-\xi-t/V_k(\xi)}(k)-\varphi_{%-\xi-t/V_k(\xi)}(k/c)}{\varphi_{-\xi-t/V_k(\xi)}(k/c)} \right ]. \]
%Furthermore, since for $n$ large enough and for all $\xi \in \R$,
%\[\frac{\varphi_{-\xi-t/V_k(\xi)}(k)-\varphi_{-\xi-t/V_k(\xi)}(k/c)}{\varphi_{%-\xi-t/V_k(\xi)}(k/c)} > 0, \]
%and since $\varphi_{\xi}(1/c)<0$ for all $\xi \in \R$, we have,
 Routine calculations yield:
\[ F_n(t)=\proba \left [ - \frac{\varphi_{\xi}(1/c)}{Z_n} \leq t_n \right ], \]
with
\[ t_n=-\varphi_{\xi}(1/c) \frac{\varphi_{-\xi-t/V_k(\xi)}(k/c)}{\varphi_{-\xi-t/V_k(\xi)}(k)-\varphi_{-\xi-t/V_k(\xi)}(k/c)}. \]
Remarking that $\varphi_{\xi}^{*}$ is an increasing function for $\xi \geq 0$ and decreasing for $\xi<0$, we have,
\begin{equation}
\label{fnt}
F_n(t) = \left \{
\begin{array}{l l}
\proba \left [ (k')^{\delta-\ind \{ \xi=0 \}} \varphi_{\xi}^{*} \left ( - \frac {\varphi_{\xi}(1/c)}{Z_n} \right ) \leq (k')^{\delta-\ind \{ \xi=0 \}} \varphi_{\xi}^{*}(t_n) \right ] & {\rm{if}} \ \xi \geq 0, \\
\proba \left [ (k')^{\delta} \varphi_{\xi}^{*} \left ( - \frac{\varphi_{\xi}(1/c)}{Z_n} \right ) \geq (k')^{\delta} \varphi_{\xi}^{*}(t_n) \right ] & {\rm{if}} \ \xi < 0.
\end{array}
\right .
\end{equation}
The asymptotic behavior of the left hand side random term 
\[ (k')^{\delta-\ind \{ \xi=0 \}} \varphi_{\xi}^{*} \left ( - \frac
{\varphi_{\xi}(1/c)}{Z_n} \right ) \]
is given by Lemma \ref{lemcvloiZ}. Let us now focus on the right hand
side deterministic term $(k')^{\delta-\ind \{\xi =0\}} \varphi_{\xi}^{*}(t_n)$. Different cases have to be considered: \\
{\underline{If $\xi >0$}}, the following sequence of asymptotic
equivalences holds
\[ t_n \sim k^{\xi+t/V_k(\xi)} \frac{\varphi_{\xi}(1/c)}{1-c^{\xi+t/V_k(\xi)}} \sim k^{\xi+t/V_k(\xi)} c^{-\xi}/\xi. \]
Since $V_{k}(\xi) \sim \ln(k)/\xi$, we have that $k^{t/V_k(\xi)} \to
{\rm{e}}^{t \xi}$ as $n \to \infty$. Thus, $t_n \sim (k')^{\xi}/\xi{\rm{e}}^{t \xi} \to \infty$ as $n \to \infty$. Lemma \ref{lemphi} III) implies that $\varphi_{\xi}^{*}(t_n) \sim 1+(k')^{\xi}{\rm{e}}^{t \xi} \sim (k')^{\xi}{\rm{e}}^{t \xi}$. Thus,
\begin{equation}
\label{casA}
\lim_{n \to \infty} (k')^{\delta-\ind \{\xi =0\}}\varphi_{\xi}^{*}(t_n) = {\rm{e}}^{t \xi}.
\end{equation}
{\underline{If $\xi =0$}}, we have
\[ t_n=\ln(c) \frac{c^{t/\ln^2(k)}-k^{t/\ln^2(k)}}{1-c^{t/\ln^2(k)}}. \]
Using the expansions,
\[ c^{t/\ln^2(k)}=1+t\frac{\ln(c)}{\ln^2(k)} + O \left ( \frac{1}{\ln^4(k)} \right ) \ {\rm{and}} \ k^{t/\ln^2(k)}=1+t\frac{1}{\ln(k)} +t^2\frac{1}{2\ln^2(k)}+o \left ( \frac{1}{\ln^2(k)} \right ), \]
we find that
\[ t_n = \ln(k) \left [ 1 + \frac{t/2-\ln(c)}{\ln(k)} + o \left ( \frac{1}{\ln(k)} \right ) \right ]. \]
Since $t_n \to \infty$ and $\ln(k)-t_n \to -t/2+\ln(c)$ as $n \to \infty$, Lemma~\ref{lemphi} II) implies that
\[ \varphi_{\xi}^{*}(t_n) \sim \varphi_{\xi}^{*}(\ln(k)) \exp[t/2-\ln(c)] = k'\exp(t/2). \]
Thus,
\begin{equation}
\label{casB}
\lim_{n \to \infty} (k')^{\delta-\ind \{\xi =0\}}\varphi_{\xi}^{*}(t_n) = {\rm{e}}^{t/2}.
\end{equation}
When {\underline{$\xi <0$}}, we have:
\[ t_n=-\frac{1}{\xi} \frac{c^{-\xi}-1}{c^{-\xi}-c^{t/V_k(\xi)}}
c^{t/V_k(\xi)}[1-(k')^{\xi+t/V_k(\xi)}]. \]
Remarking that
\[ c^{t/V_k(\xi)}=1+t\frac{\ln(c)}{V_k(\xi)} + O \left (
\frac{1}{V_k(\xi)} \right ), \]
and that $(k')^{t/V_k(\xi)}=1+o(1)$ lead to the following expansion:
\begin{eqnarray*}
t_n & = & -\frac{1}{\xi} \left [ 1 + t\frac{\ln(c)}{V_k(\xi)} \right ] \left [ 1+o\left( \frac{1}{V_k(\xi)} \right ) \right ] \left [ 1 + t \frac{\ln(c)}{(c^{-\xi}-1)V_k(\xi)} + o\left( \frac{1}{V_k(\xi)} \right ) \right ] \left [ 1-(k')^{\xi}+o(k^{\xi}) \right ] \\
 \ & = & u_n(1+\varepsilon_n),
 \end{eqnarray*}
with $u_n=-1/\xi[1+t\ln(c)/V_k(\xi)]$ and
\[ \varepsilon_n = t \frac{\delta \ln(c)}{c^{-\xi}-1}k^{-\delta} - (k')^{\xi} + o(k^{-\delta}) + o(k^{\xi}). \]
Two situations have to be considered: \\
\noindent {\underline{If $-1/2 \leq \xi <0$}}, it follows that
\[ \lim_{n \to \infty} \frac{\varepsilon_n}{\varphi_{\xi}^{*}(u_n)} = - \frac{1}{t\xi\ln(c)} \left [ \frac{t\xi\ln(c)}{c^{-\xi}-1}+c^{-\xi} \right ], \]
and Lemma \ref{lemphi} IV) i) implies
\[ \varphi_{\xi}^{*}(t_n) \sim \varphi_{\xi}^{*}(u_n) \left [1+\frac{1}{t\xi\ln(c)} \left ( \frac{t\xi\ln(c)}{c^{-\xi}-1}+c^{-\xi} \right ) \right ] \]
and thus
\begin{equation}
\label{casC}
\lim_{n \to \infty} (k')^{\delta-\ind \{\xi =0\}}\varphi_{\xi}^{*}(t_n) = 1+ t \frac{\xi \ln(c)}{c^{-\xi}-1}.
\end{equation}
{\underline{If $\xi<-1/2$}}, $\varepsilon_n/\varphi_{\xi}^{*}(u_n) \to (1-c^{-\xi})^{-1}$ as $n \to \infty$. Thus, Lemma \ref{lemphi} IV) i) yields
\begin{equation}
\label{casD}
\lim_{n \to \infty} (k')^{\delta-\ind \{\xi =0\}}\varphi_{\xi}^{*}(t_n) = -\frac{t}{2} \ln(c) \frac{c^{-\xi-1/2}}{1-c^{-\xi}}.
\end{equation}
Collecting (\ref{fnt})-(\ref{casD}) with Lemma \ref{lemcvloiZ}
concludes the proof. \findemo 

\section{Simulation study}
\label{secsimul}

In this section, the improvement brought by the bias correction is illustrated through a simulation study. Next, a comparison with
classical extreme value estimators is proposed. For each of the
distributions considered in this section, $N=100$ random samples of
size $n=500$ are generated.  

\subsection{Bias corrected estimator behavior}
\label{subsecsimul1}

\begin{figure}[h]
\begin{center}
\begin{tabular}{cc}
\subfigure[]{\psfig{figure=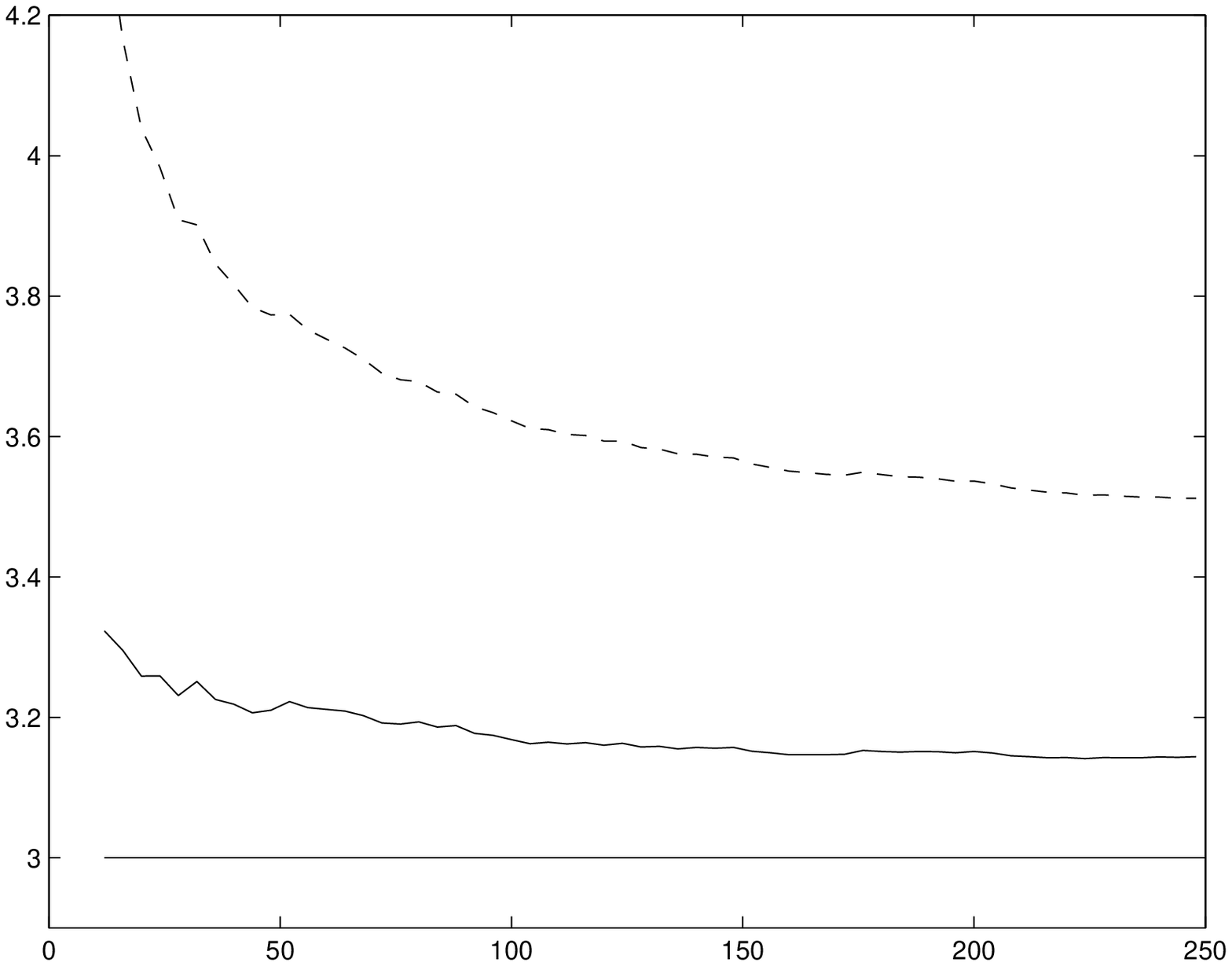,width=6.5cm}} &
\subfigure[]{\psfig{figure=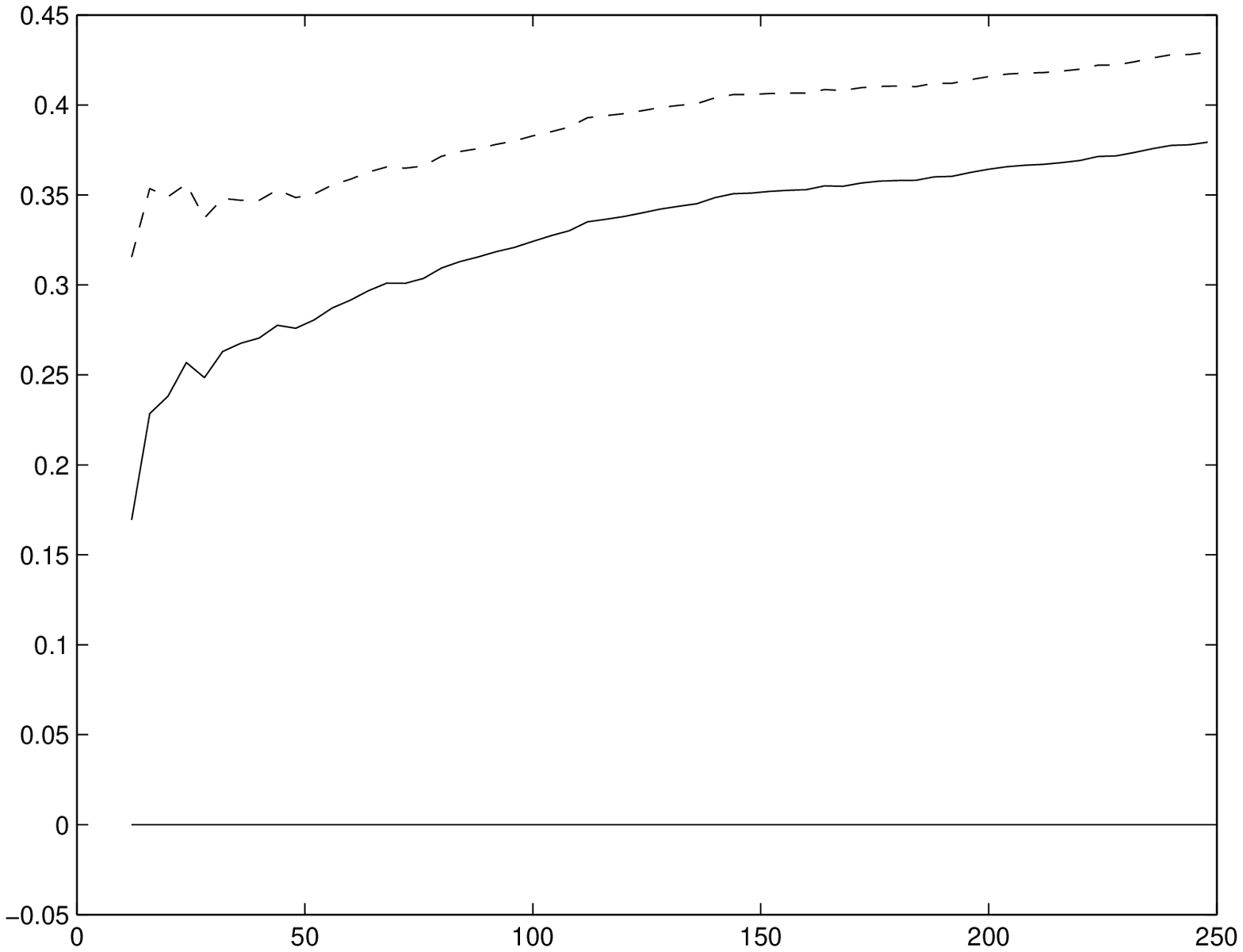,width=6.5cm}} \\
\subfigure[]{\psfig{figure=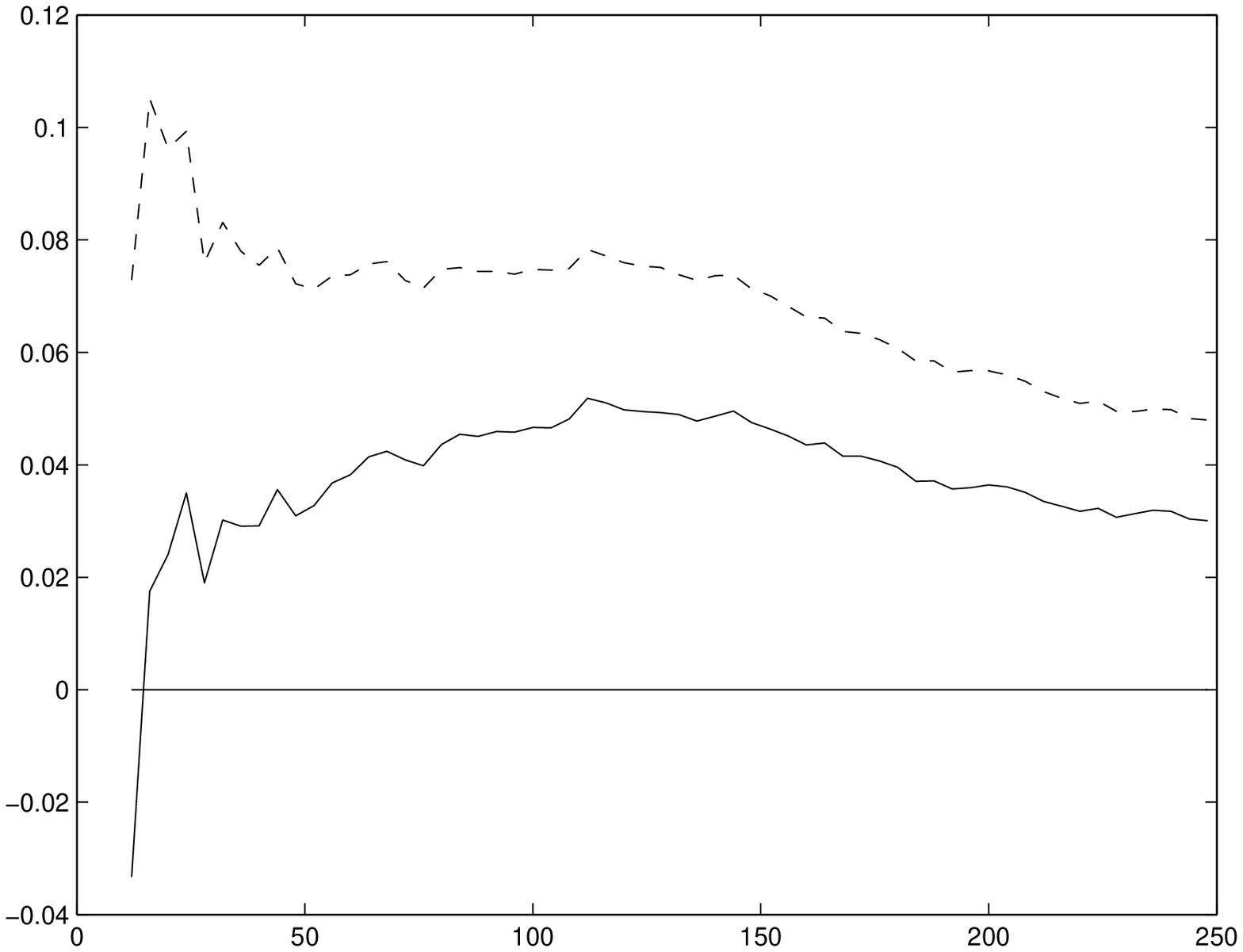,width=6.5cm}} &
\subfigure[]{\psfig{figure=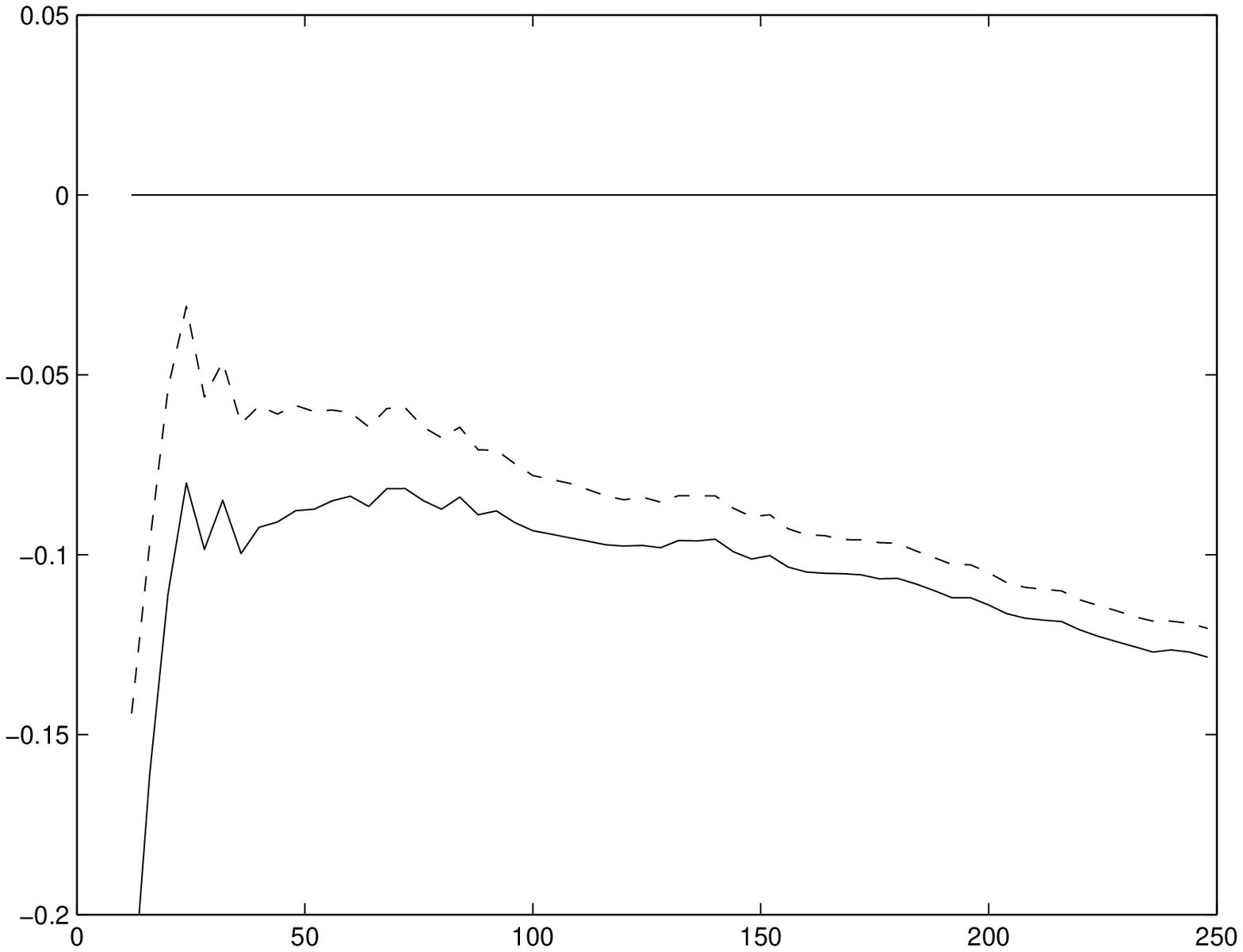,width=6.5cm}} \\
\subfigure[]{\psfig{figure=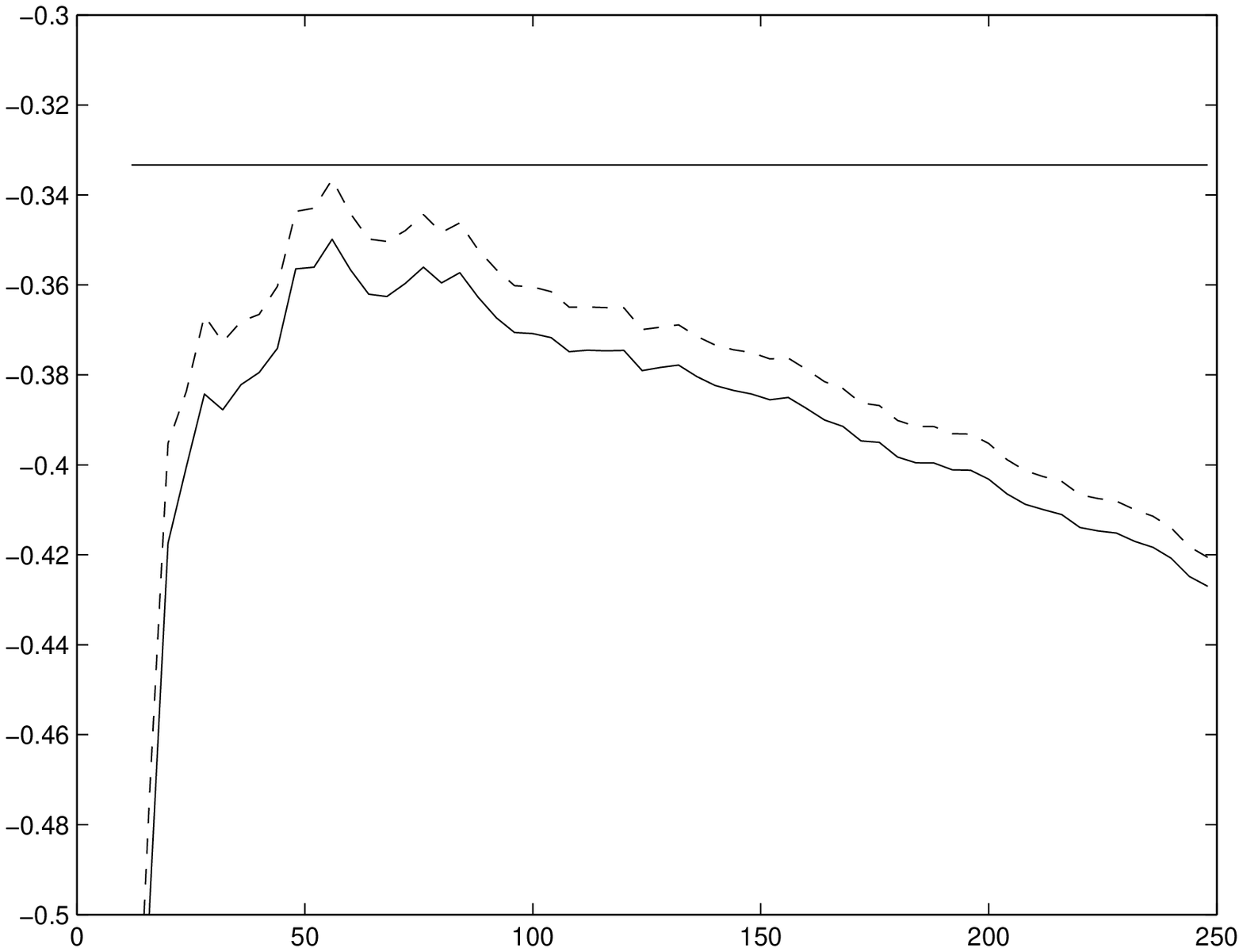,width=6.5cm}} & 
\subfigure[]{\psfig{figure=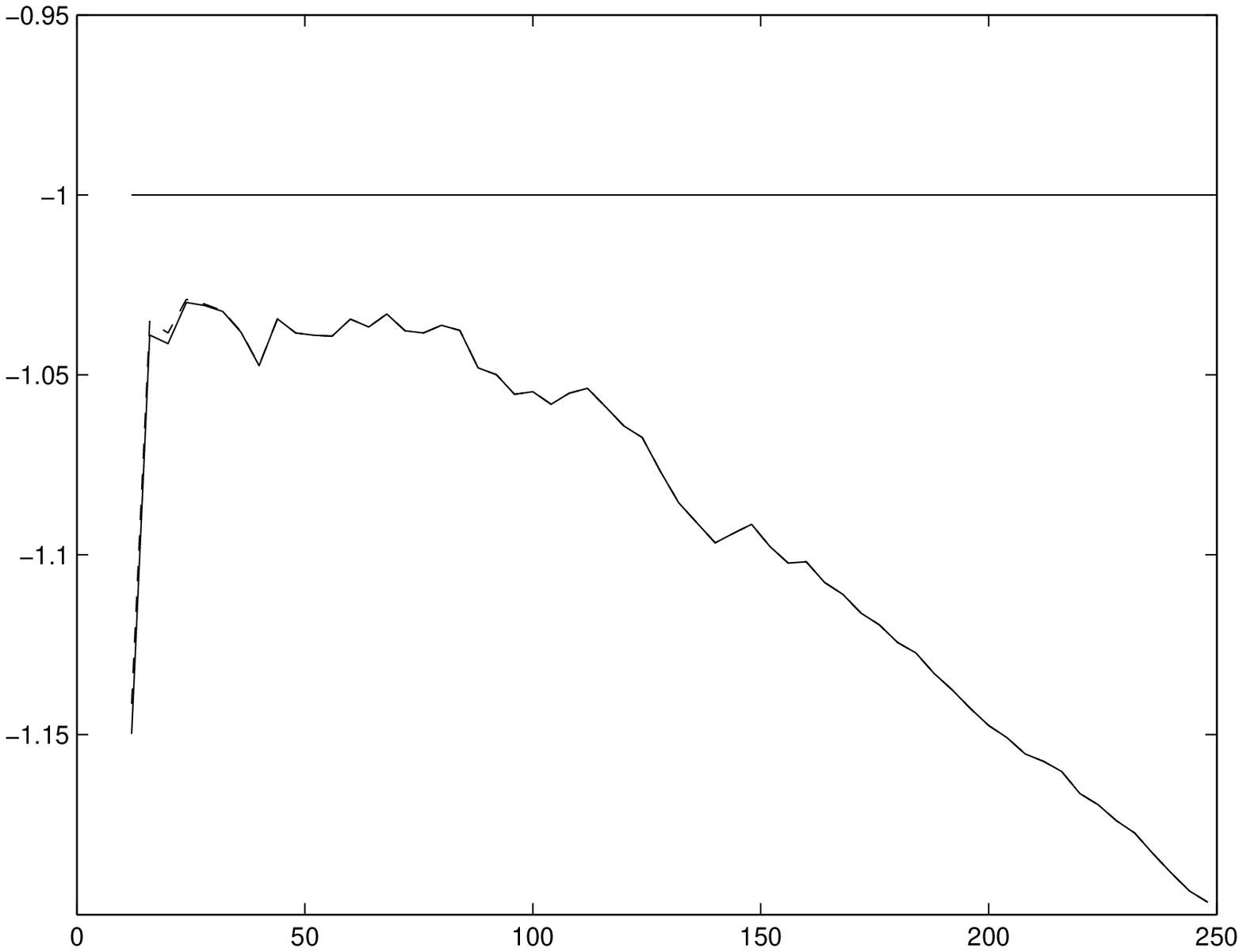,width=6.5cm}} \\
\end{tabular}
\end{center}
\caption{Comparison of the bias corrected estimator
${\hat{\xi}}_{k,n}^*$ (solid line) and ${\hat{\xi}}_{k,n}$ (dashed
line) for (a) the Fr\'{e}chet distribution, (b) the Weibull distribution
with $(\lambda,\tau)=(1,1/2)$, (c) the Weibull distribution
with $(\lambda,\tau)=(1,1)$, (d) the Weibull distribution
with $(\lambda,\tau)=(1,3/2)$, (e) the Weibull$_{{\rm{M}}}$
distribution with $\xi=-1/3$ and (f) the Weibull$_{{\rm{M}}}$
distribution with $\xi=-1$.}
\label{fig1}
\end{figure}

We first study the behavior of the bias corrected estimator
${\hat{\xi}}_{k,n}^*$ versus the estimator ${\hat{\xi}}_{k,n}$. In
this aim, the following distributions are considered: (see Table
\ref{tab1} for their parameterizations)
\begin{itemize}
\item[$\bullet$] {\underline{Case $\xi>0$}}, Fr\'{e}chet distribution
with $\xi=3$.
\item[$\bullet$] {\underline{Case $\xi=0$}}, Weibull distribution with $(\lambda,\tau)=(1,1/2)$, $(\lambda,\tau)=(1,1)$ and $(\lambda,\tau)=(1,3/2)$.
\item[$\bullet$] {\underline{Case $-1/2<\xi<0$}}, Weibull$_{{\rm{M}}}$
distribution with $\xi=-1/3$.
\item[$\bullet$] {\underline{Case $\xi<-1/2$}}, Weibull$_{{\rm{M}}}$
with $\xi=-1$.
\end{itemize}
In Figure \ref{fig1}, the empirical mean over the $N$ samples of ${\hat{\xi}}_{k,n}^*$ and
${\hat{\xi}}_{k,n}$ is represented as a function of the
number $k$ of upper order statistics (``Hill plot''). The true value
of $\xi$ is represented by a straight line. To compute these estimators, we
choose $c=k/k'=4$. If $\xi > 0$ (Fr\'{e}chet distribution,
Figure~\ref{fig1} (a)), the behavior of ${\hat{\xi}}_{k,n}$ is
improved by  the bias correction. If $\xi=0$ (Weibull
distribution), the estimation is highly influenced by the parameter $\tau$
which controls the rate of convergence of the slowly varying function
$\ell$ (see table~\ref{tab1}). If $\tau \leq 1$ (Figure~\ref{fig1}
(b), (c)), ${\hat{\xi}}_{k,n}^*$ is less biased than
${\hat{\xi}}_{k,n}$ and if $\tau>1$ (Figure~\ref{fig1}
(d)), the bias correction does not improve the behavior of
${\hat{\xi}}_{k,n}$. If $-1/2<\xi<0$ (Figure~\ref{fig1} (e)),
${\hat{\xi}}_{k,n}^*$ is slightly more biased than the estimator
${\hat{\xi}}_{k,n}$. Finally, if $\xi<-1/2$ (Weibull$_{{\rm{M}}}$
distribution, Figure~\ref{fig1} (f)), there is no correction (${\hat{\xi}}_{k,n}={\hat{\xi}}_{k,n}^*$). As a conclusion, it seems that the bias
correction improves (or at least does not really degrade) the
behavior of our estimator. Thus, in the sequel, we
focus on the behavior of ${\hat{\xi}}_{k,n}^*$. 

\subsection{Comparison with other estimators}

\begin{figure}[h]
\begin{center}
\begin{tabular}{c c}
\subfigure[mean]{\psfig{figure=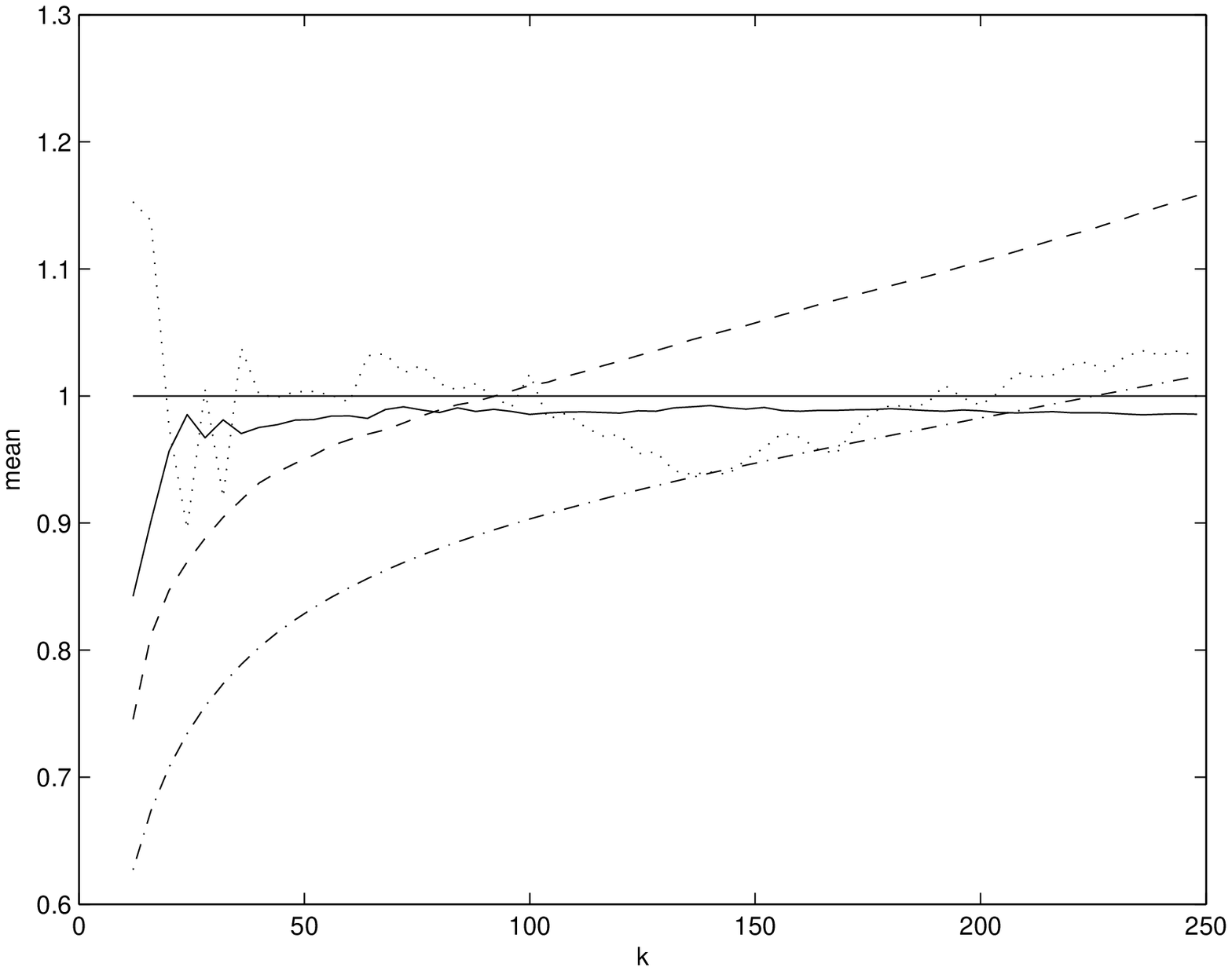,width=6.5cm}} &
\subfigure[MSE]{\psfig{figure=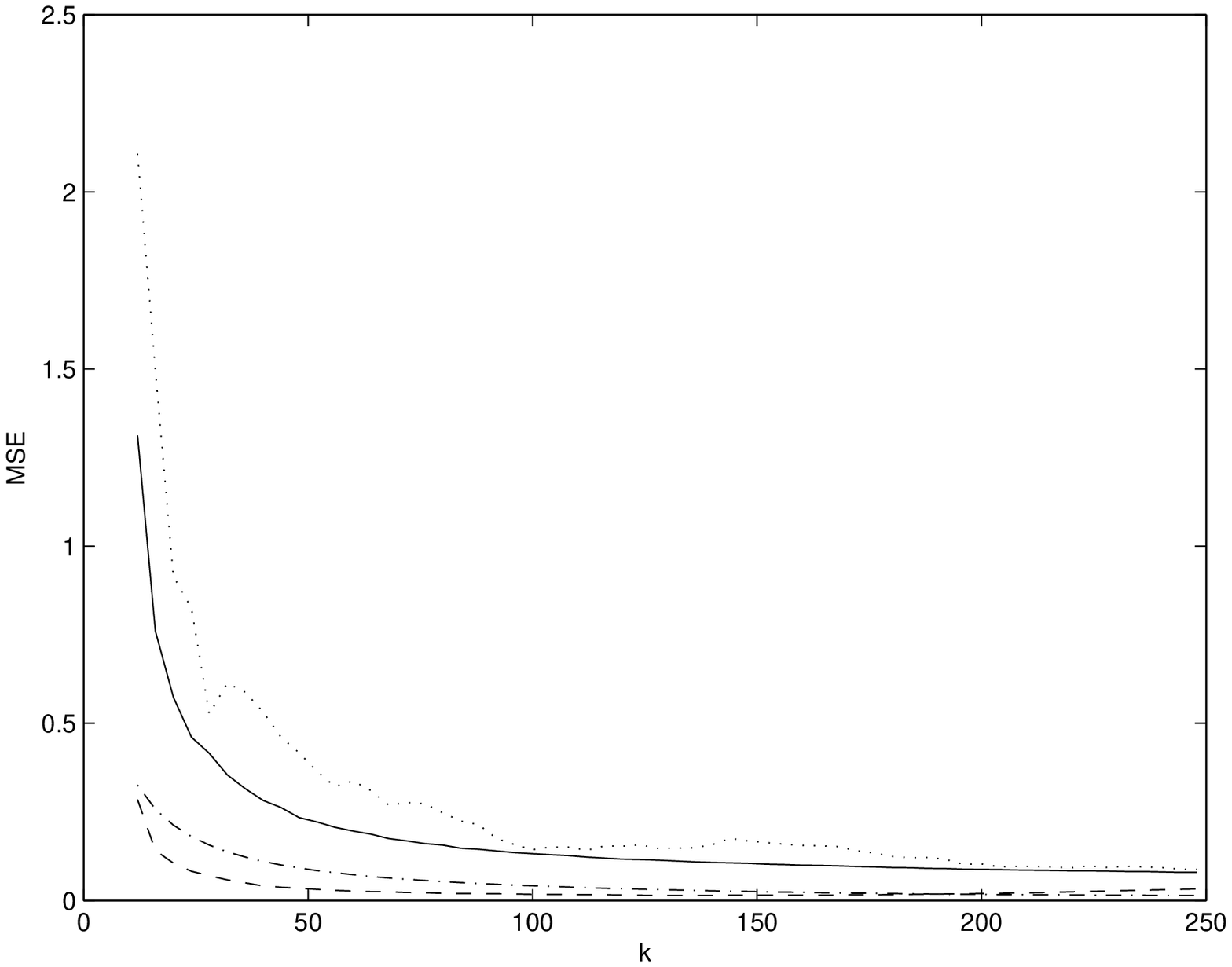,width=6.5cm}} \\
\subfigure[mean]{\psfig{figure=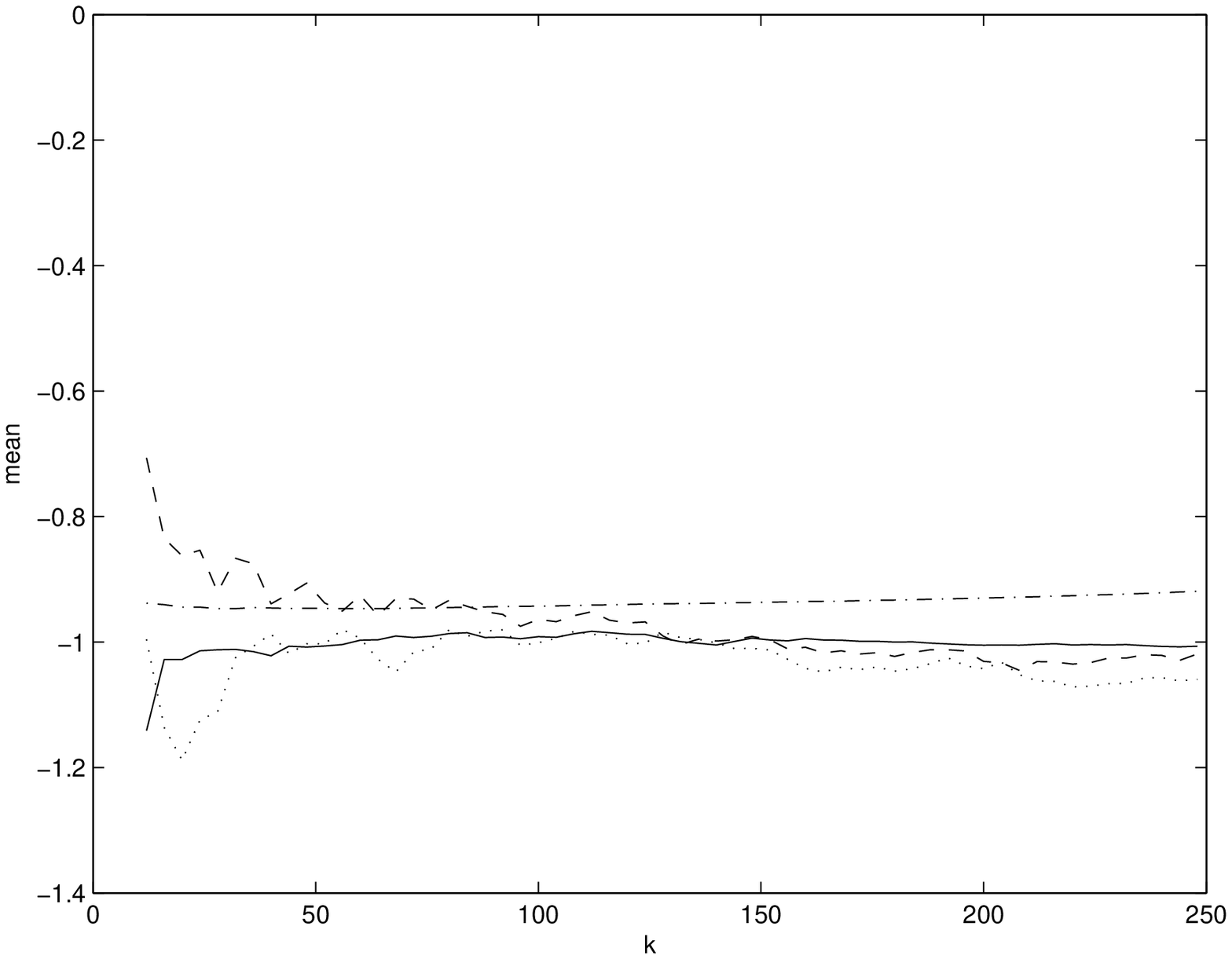,width=6.5cm}} &
\subfigure[MSE]{\psfig{figure=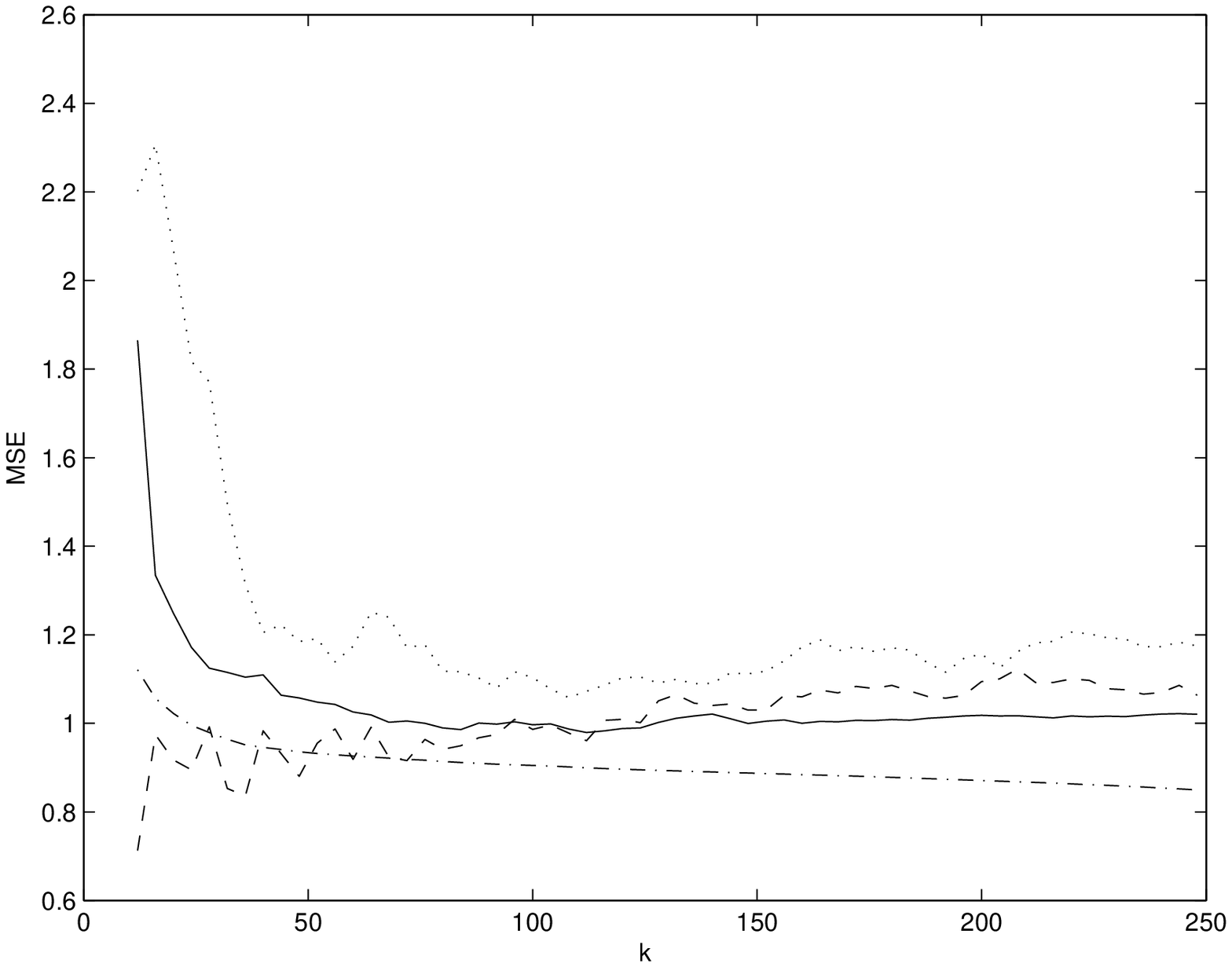,width=6.5cm}} \\
\subfigure[mean]{\psfig{figure=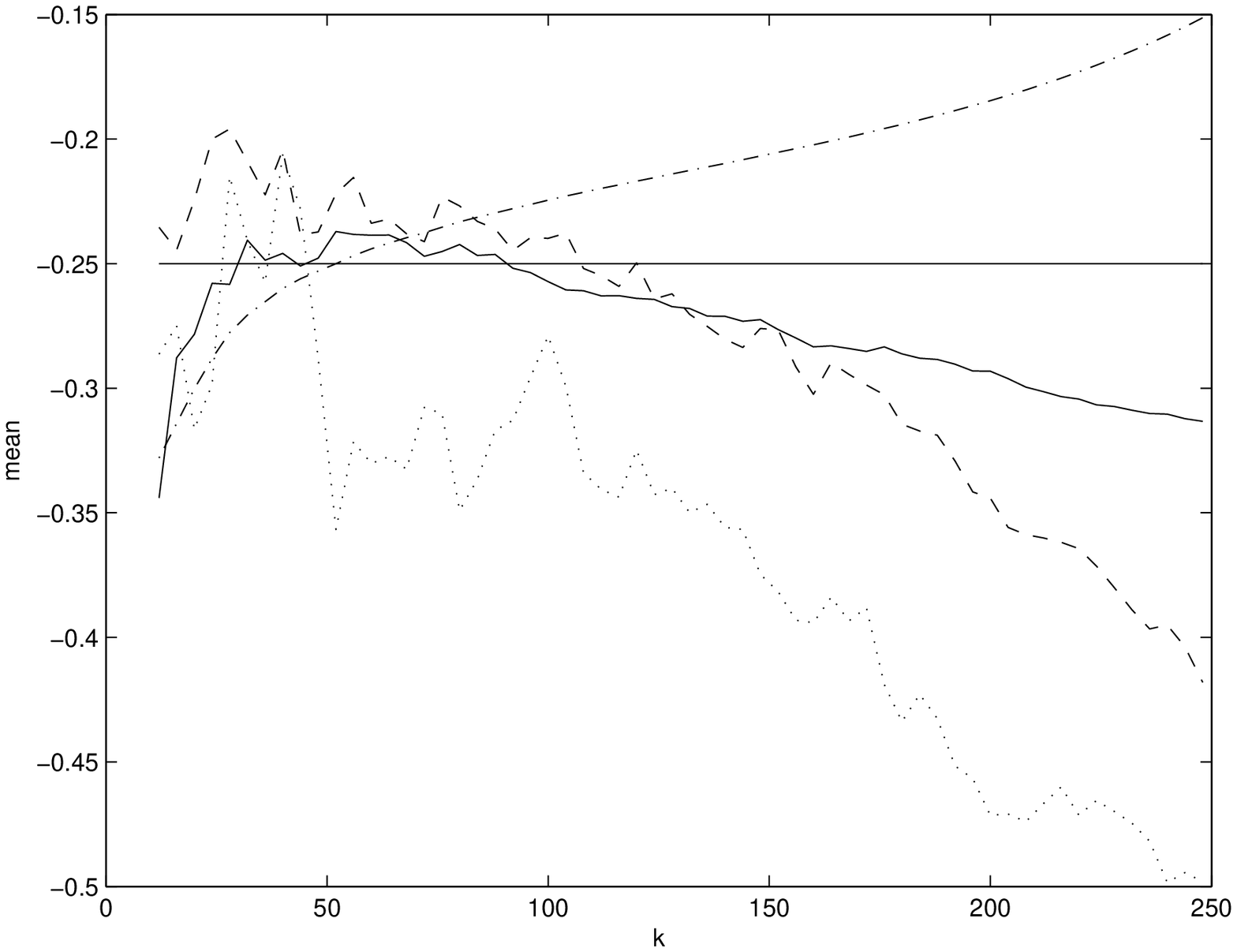,width=6.5cm}} &
\subfigure[MSE]{\psfig{figure=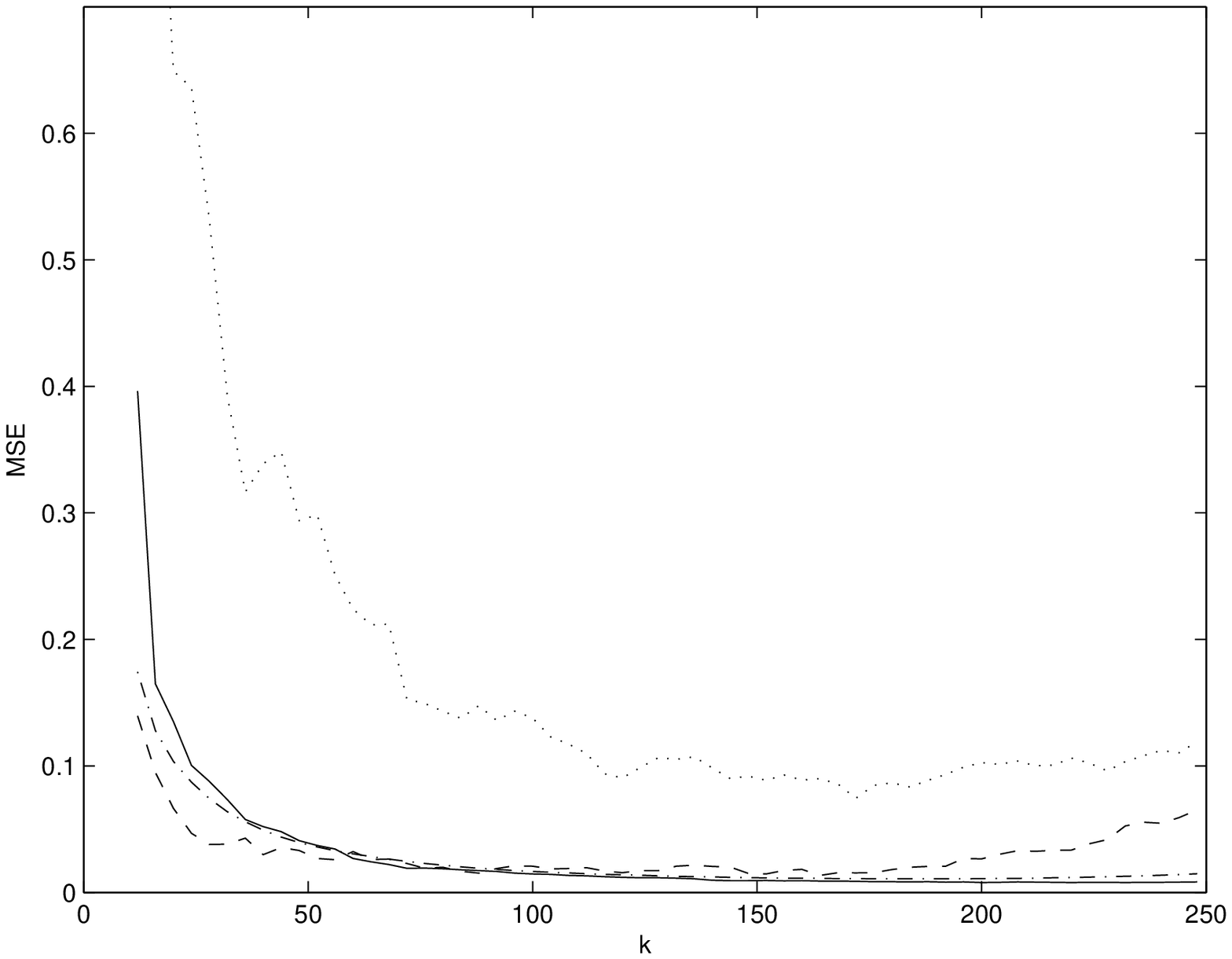,width=6.5cm}} \\
\end{tabular}
\end{center}
\caption{Comparison of the estimator ${\hat{\xi}}_{k,n}^*$ (solid
line), moment estimator (dashed line), Pickands'
estimator (dotted line) and the generalized Zipf estimator
(dash-dot line) for (a), (b) the Burr distribution with
$(\beta,\tau,\lambda)=(1,1,1)$, (c), (d) the standard normal
distribution and (e), (f) the Weibull$_{{\rm{M}}}$ distribution with
$\xi=-1/4$.}
\label{fig2}
\end{figure}

\begin{figure}
\begin{center}
\begin{tabular}{c c}
\subfigure[Mean]{\psfig{figure=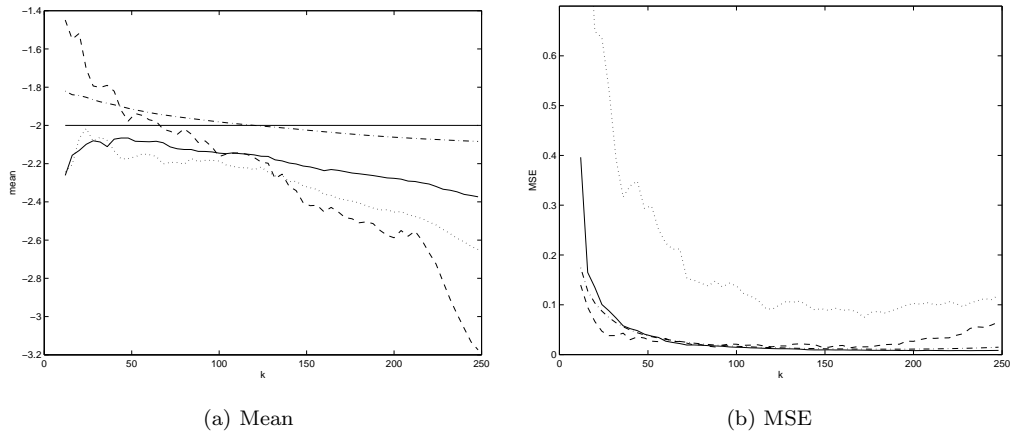,width=6.5cm}} &
\subfigure[MSE]{\psfig{figure=MSECCE3.eps,width=6.5cm}} \\
\end{tabular}
\end{center}
\caption{Comparison of the estimator ${\hat{\xi}}_{k,n}^*$ (solid
line), moment estimator (dashed line), Pickands'
estimator (dotted line) and the generalized Zipf estimator
(dash-dot line) for the Weibull$_{{\rm{M}}}$
distribution with $\xi=-2$.}
\label{fig3}
\end{figure}

\begin{figure}
\begin{center}
\begin{tabular}{c c}
\subfigure[mean]{\psfig{figure=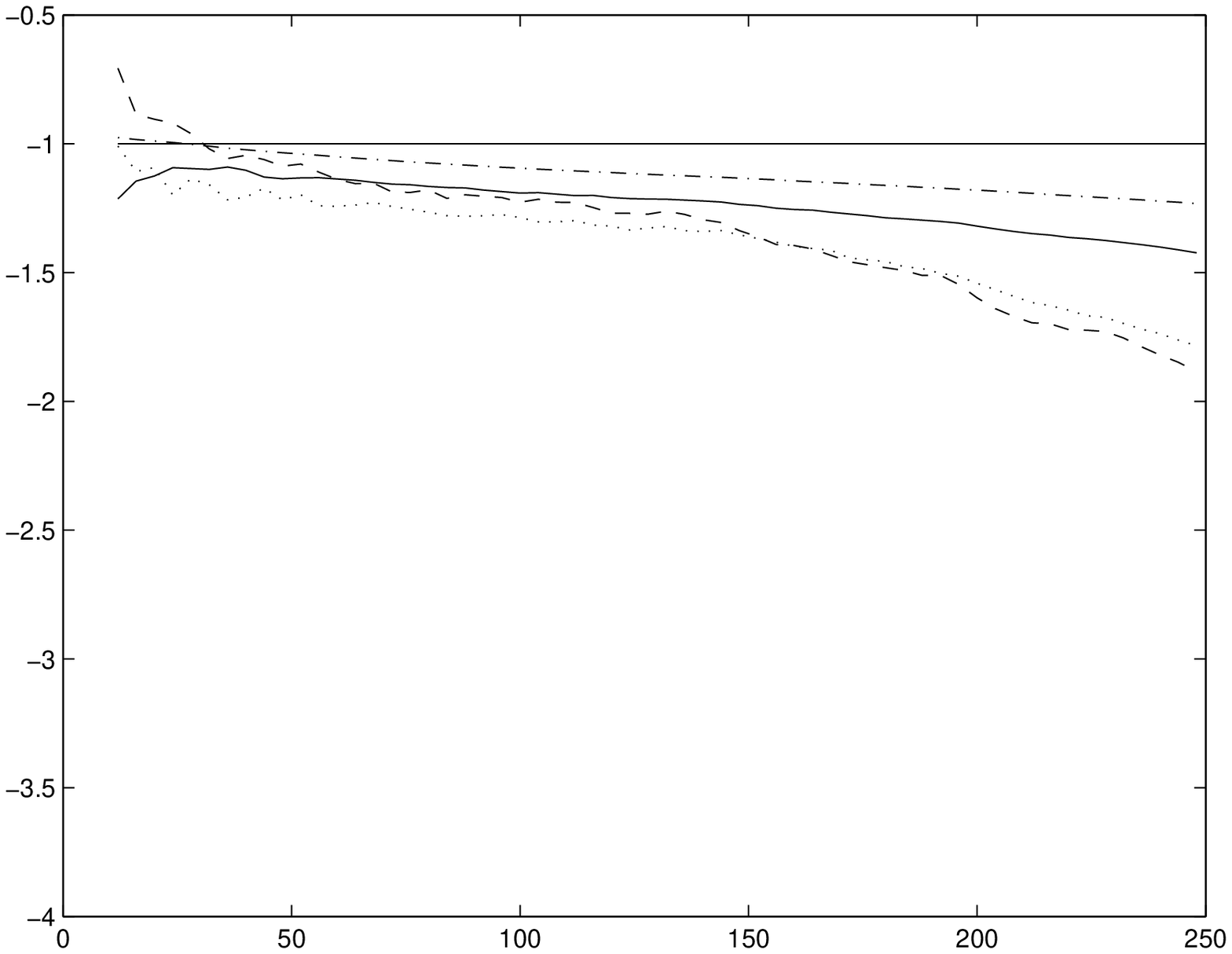,width=6.5cm}} &
\subfigure[MSE]{\psfig{figure=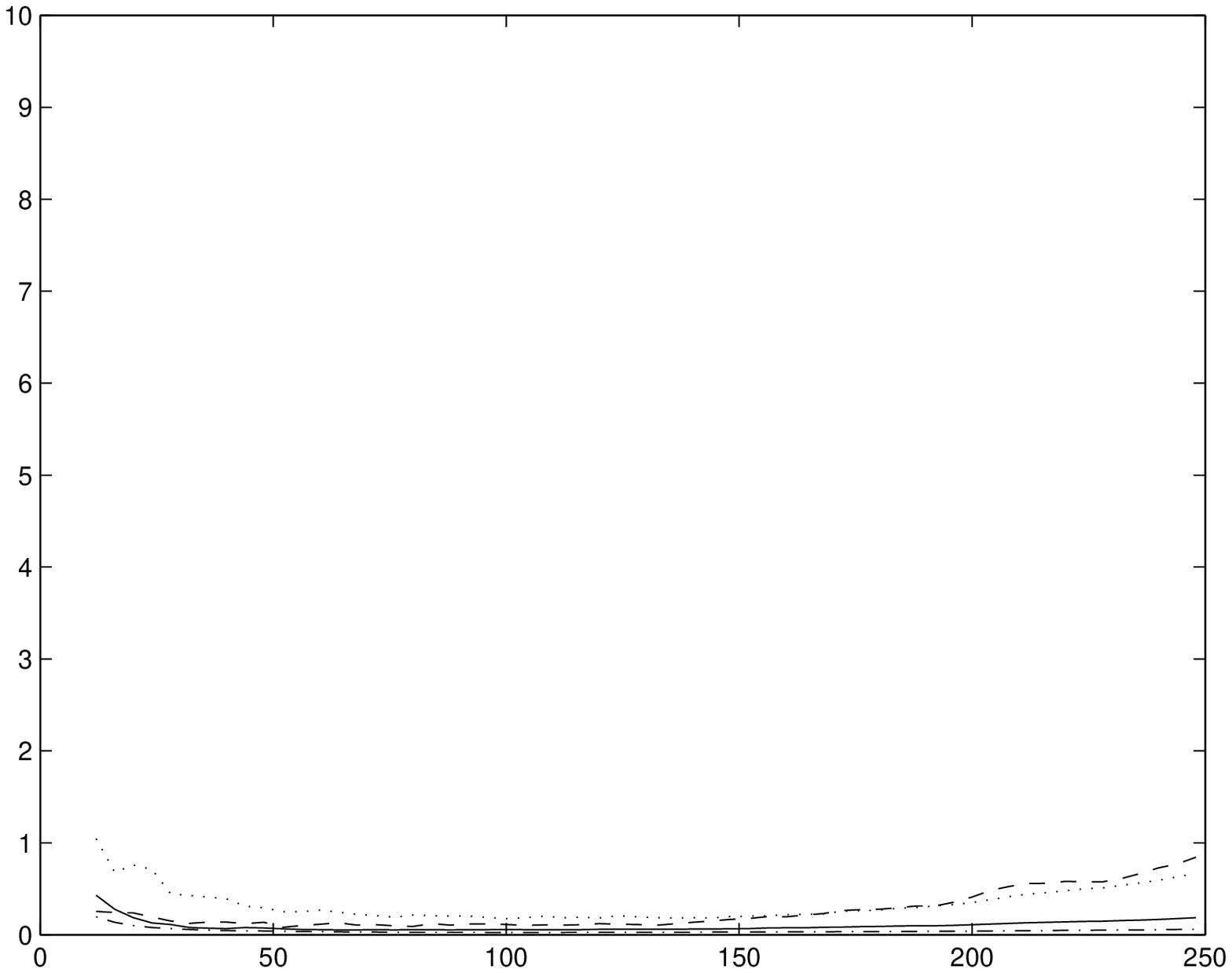,width=6.5cm}} \\
\subfigure[mean]{\psfig{figure=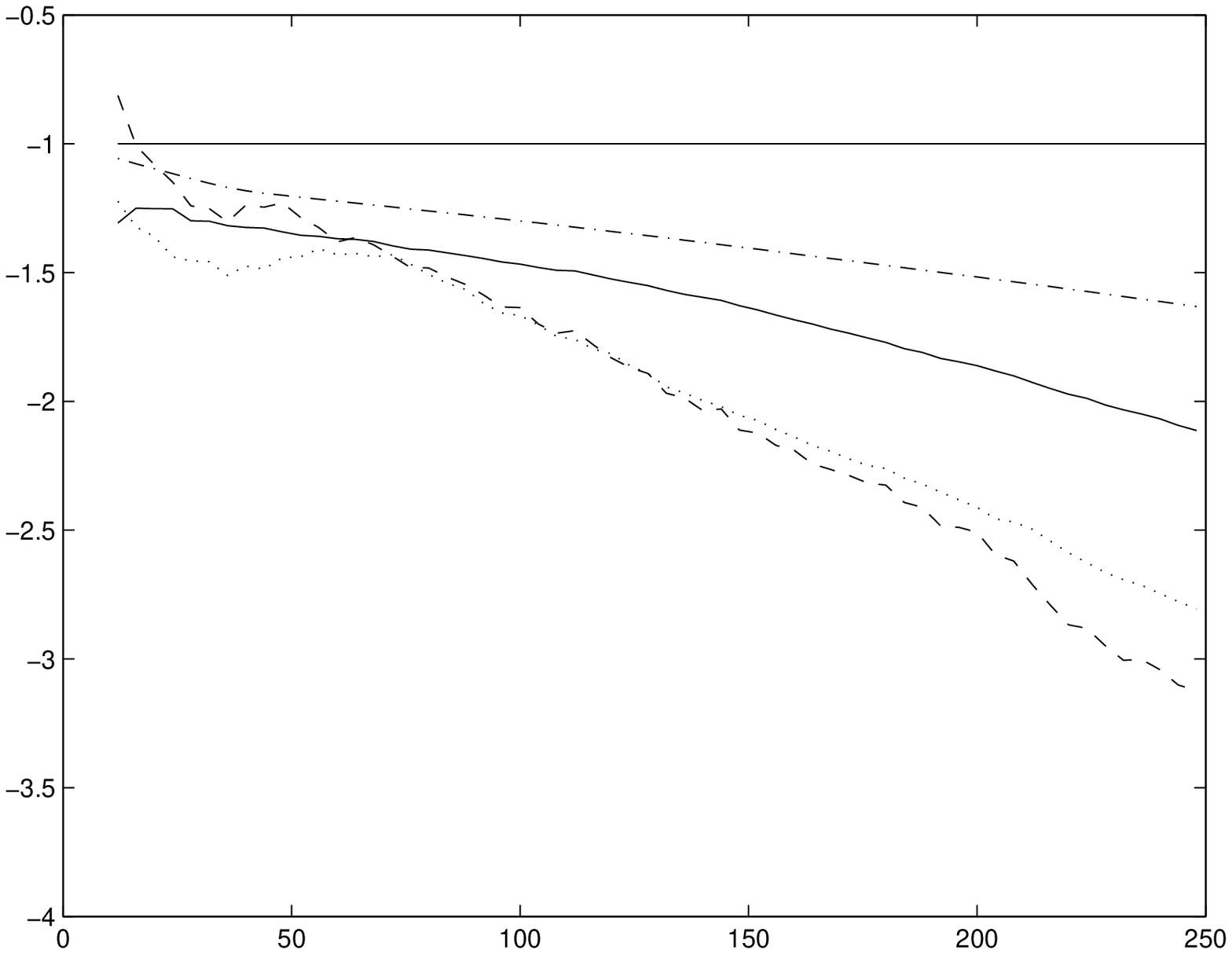,width=6.5cm}} &
\subfigure[MSE]{\psfig{figure=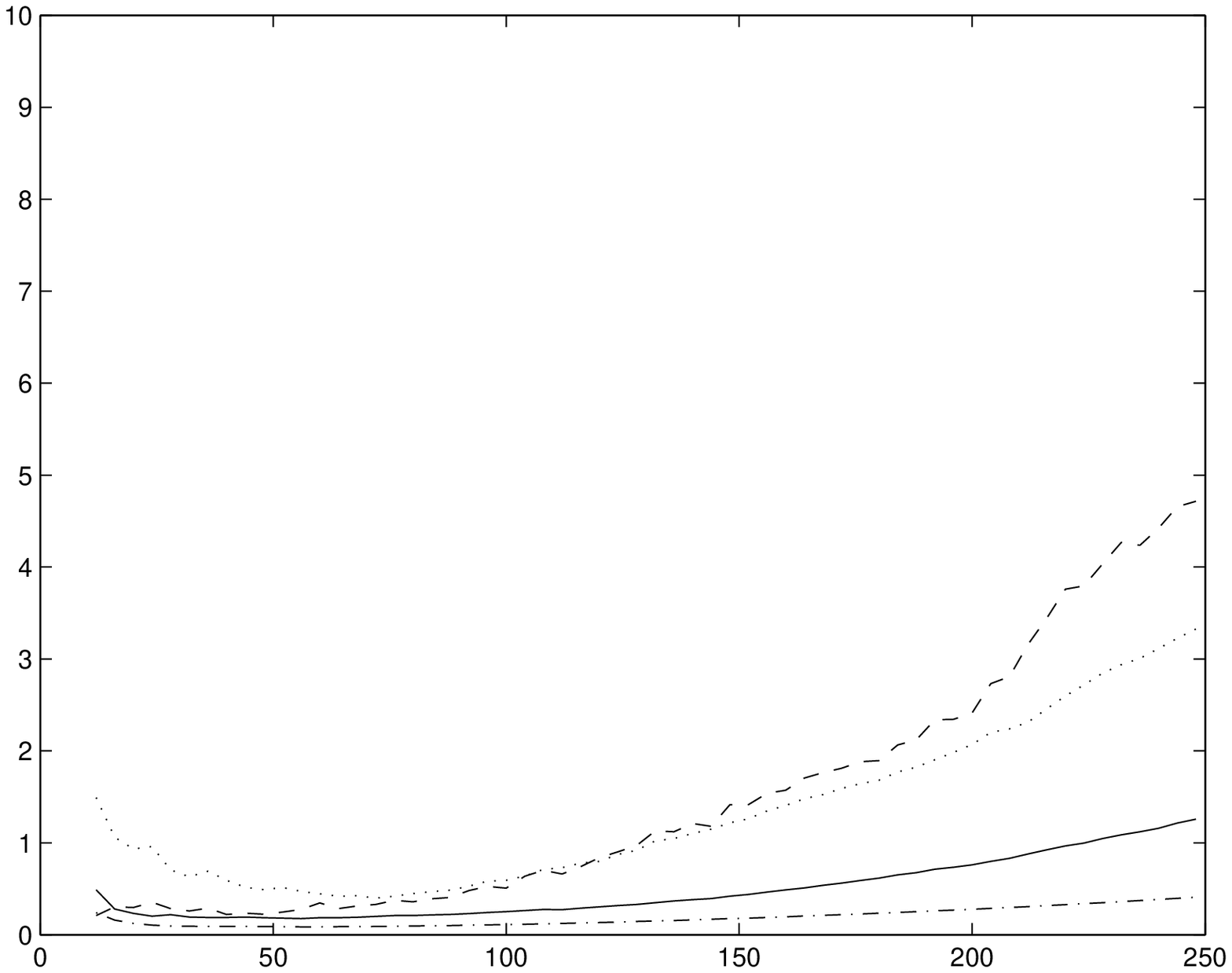,width=6.5cm}} \\
\subfigure[mean]{\psfig{figure=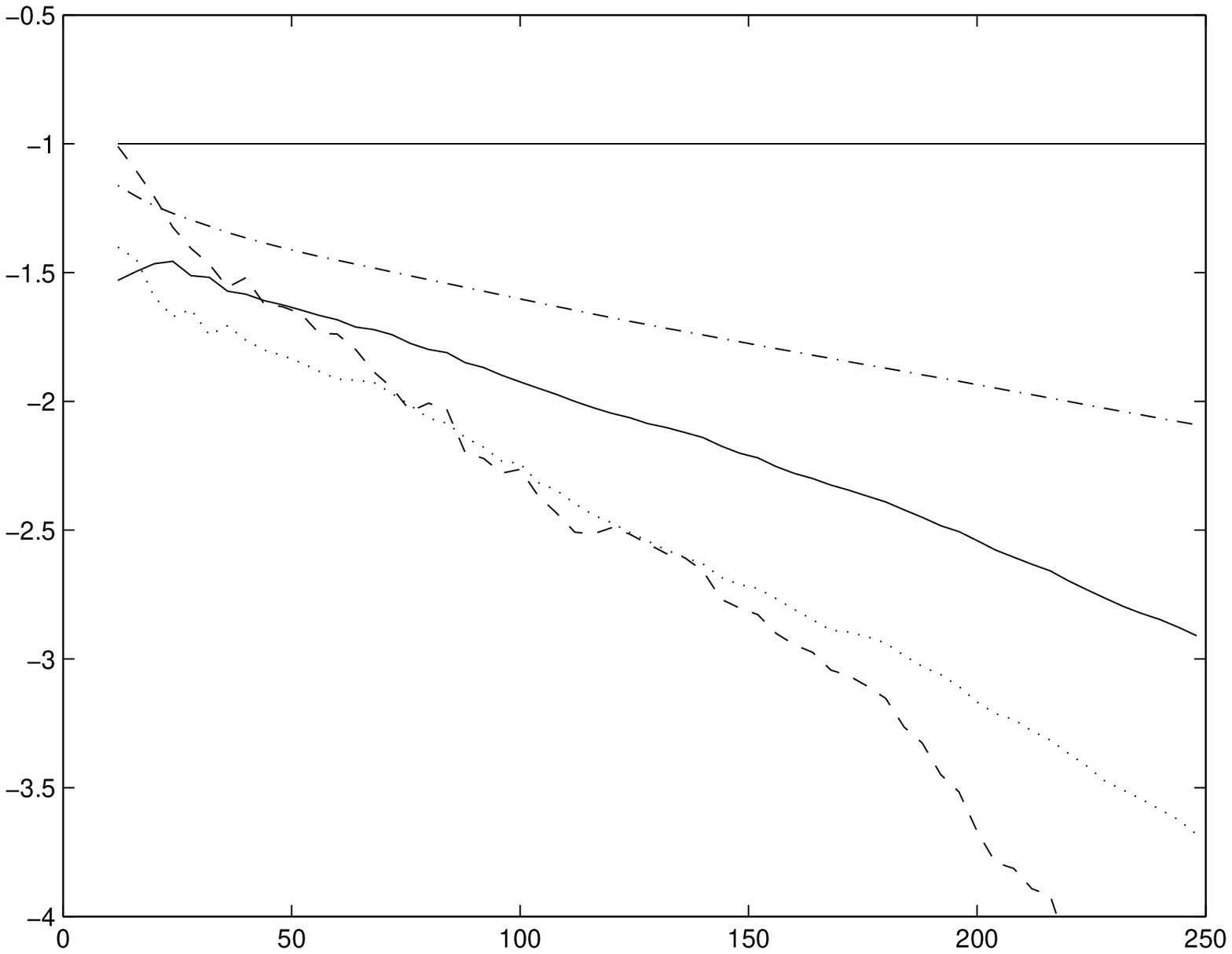,width=6.5cm}} &
\subfigure[MSE]{\psfig{figure=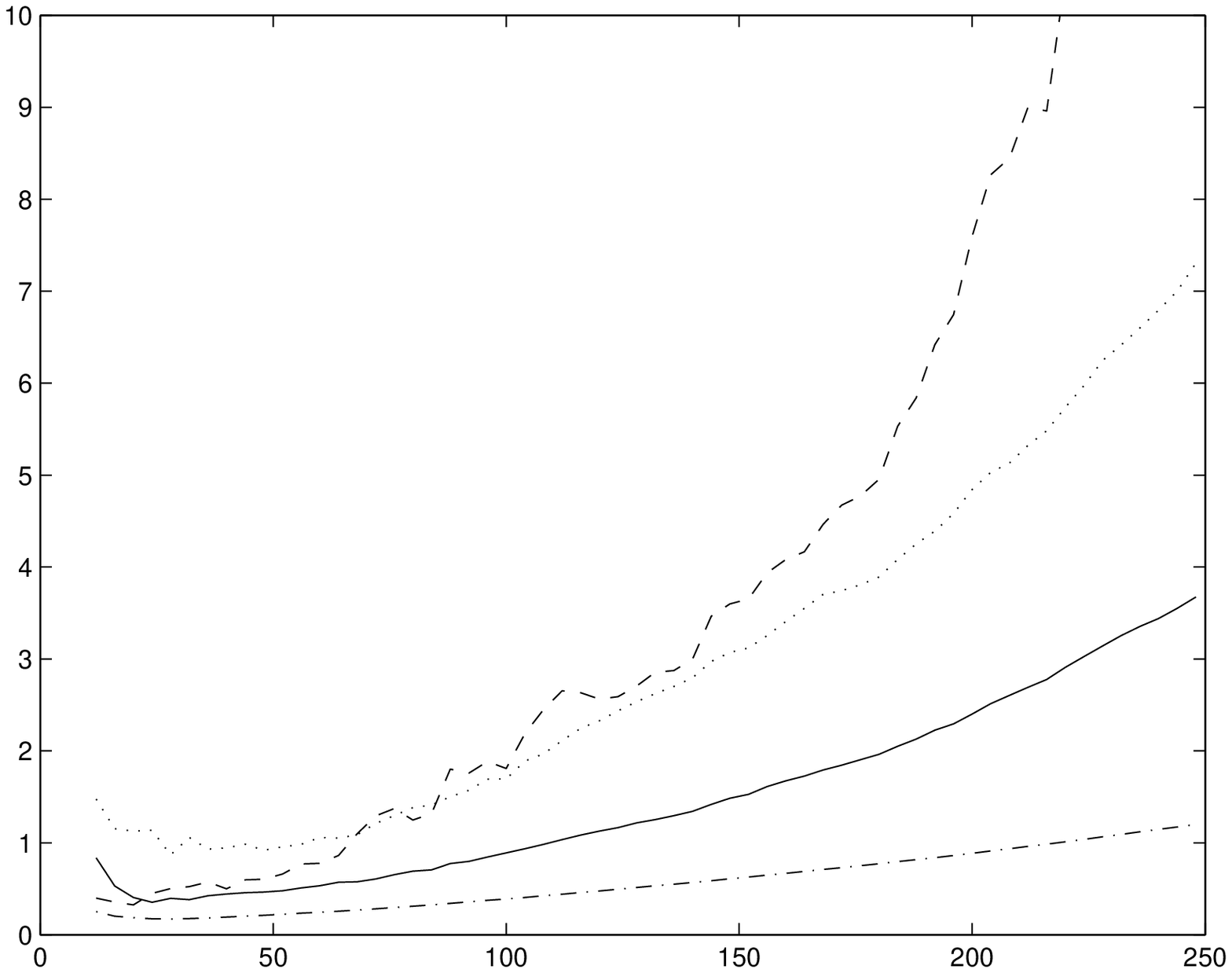,width=6.5cm}} \\
\end{tabular}
\end{center}
\caption{Comparison of the estimator ${\hat{\xi}}_{k,n}^*$ (solid
line), moment estimator (dashed line), Pickands'
estimator (dotted line) and the generalized Zipf estimator
(dash-dot line) for the reversed Burr distribution with (a), (b)
$\lambda=1$, (c), (d) $\lambda=2$ and (e), (f) $\lambda=3$.}
\label{fig4}
\end{figure}

The estimator ${\hat{\xi}}_{k,n}^*$ is now compared with the following
well known estimators: {\bf{Pickands' estimator}} ${\hat{\xi}}_{k,n}^P$,
the {\bf{moment estimator}} proposed by Dekkers, Einmahl and de Haan
\cite{dekeinhaa89} and defined by:
\[ {\hat{\xi}}_{k,n}^{M} = {\hat{\xi}}_{k,n}^{H} + 1 - \frac{1}{2}  \left [ 1
- \frac{({\hat{\xi}}_{k,n}^{H})^2}{S_{k,n}} \right ]^{-1}, \]
where $S_{k,n} = 1/k \sum_{i=1}^{k}
[\ln(X_{n-i+1,n})-\ln(X_{n-k,n})]^2$ and the {\bf{generalized Zipf estimator}}
\cite{beidiegui03} defined by:
\[ {\hat{\xi}}_{k,n}^Z= \frac{\sum_{j=1}^k \ln[ (k+1)/j] \ln UH_{j,n}
-\frac{1}{k} \sum_{j=1}^k \ln[ (k+1)/j] \sum_{j=1}^k \ln
UH_{j,n}}{\sum_{j=1}^k \ln^2[(k+1)/j] - \frac{1}{k} \left (
\sum_{j=1}^k \ln [(k+1)/j] \right )^2}, \]
with
\[ UH_{j,n}=X_{n-j,n} \left ( \frac{1}{j} \sum_{i=1}^{j} \ln
X_{n-i+1,n} - \ln X_{n-j,n} \right ). \]
The following distributions are considered: (see table \ref{tab1} for
their parameterization)
\begin{itemize}
\item[$\bullet$] {\underline{Case $\xi>0$}}, Burr distribution for which $\xi=1/(\tau \lambda)$ with $(\beta,\tau,\lambda)=(1,1,1)$.
\item[$\bullet$] {\underline{Case $\xi=0$}}, standard normal distribution.
\item[$\bullet$] {\underline{Case $-1/2<\xi<0$}}, Weibull$_{{\rm{M}}}$
distribution with $\xi=-1/4$.
\item[$\bullet$] {\underline{Case $\xi<-1/2$}}, Weibull$_{{\rm{M}}}$
with $\xi=-2$.
\end{itemize}
In Figures \ref{fig2}-\ref{fig4}, the empirical mean and the empirical Mean Squared Error
(MSE) of each estimator are represented as functions of $k$ and we
also choose $c=4$. If $\xi >0$ (Burr distribution),
${\hat{\xi}}_{k,n}^*$ is less biased than the other estimators
(Figure~\ref{fig2} (a)) but it suffers from a high variance (Figure~\ref{fig2} (b)). If $\xi=0$ (Gaussian distribution, Figure~\ref{fig2}
(c), (d)), all estimates yield very poor results. If $-1/2<\xi<0$
(Weibull$_{\rm{M}}$ distribution, Figure~\ref{fig2} (e), (f)),
${\hat{\xi}}_{k,n}^*$ provides the best estimation and if $\xi<-1/2$
(Figure~\ref{fig3}), generalized Zipf estimator and
${\hat{\xi}}_{k,n}^*$ are equivalent from the MSE point of view. \\
\noindent Finally, let us focus on the influence of the rate of
convergence of the slowly varying function on the estimation of $\xi$.
In this aim, we consider the reversed
Burr distribution for which $\xi=-1/(\lambda \tau)$ (see Table~\ref{tab1} for its parameterization). Here, the parameter $\lambda$ controls the rate of convergence of the
slowly varying function (see Section \ref{secex}). The larger is $\lambda$, the slower
$\ell$ converges to a constant. This is illustrated in Figure~\ref{fig4} for $x_F=10$, $w=1$, $\tau=1/\lambda$ with $\lambda=1$, $\lambda=2$ and
$\lambda=3$. In all cases, ${\hat{\xi}}_{k,n}^*$
performed better than Pickands' and moment estimators and the best
estimation is provided by the generalized Zipf estimator. \\
\noindent As a conclusion, ${\hat{\xi}}_{k,n}^*$
performes well in the Weibull domain of attraction ($\xi<0$) and it is competitive with
Pickands' and moment estimator if $\xi \geq 0$.

\begin{rem}
{\rm{The high volatility of the previous Hill plots points out the importance of the
choice of the number $k$ of upper order statistics. A lot of methods
to select this parameter have been proposed (see Danielsson et
{\it{al.}} \cite{danhaapenvri97}, Gomes and Oliveira \cite{gomoli01},
Guillou and Hall \cite{guihal01}, \ldots). A part of our future work will consist in the adaptation of
these methods to our estimator. However, note that in all cases, our
estimator is less influenced by this choice than other estimators.}}
\end{rem}

\clearpage

%%%%%%%%%%%%%%%%%%%%%%%
%% DEBUT DE L'ANNEXE %%
%%%%%%%%%%%%%%%%%%%%%%%

\begin{center}
APPENDIX
\end{center}
 \ \\
\noindent {\bf{Proof of Lemma \ref{lem1}}} $-$ First, we focus on the case $\xi < 0$. Recall that one can take in relation (\ref{equiv1loimax}) $a(t)=\xi [U(t)-x_F]$ where $x_F$ is the right endpoint of the distribution function $F$. Thus,
\begin{equation}
\label{(A)}
\frac{U(tx)-U(t)}{a(t)} = -\frac{1}{\xi} - \frac{x_F-U(tx)}{\xi [U(t)-x_F]} < -\frac{1}{\xi},
\end{equation}
since $x_F>U(t)$ for all $t$. Furthermore, since $U$ is a non-decreasing function, we have that for all $A>0$ and $x \geq A$,
\[ \frac{U(tx)-U(t)}{a(t)} \geq \frac{U(tA)-U(t)}{a(t)}. \]
Using (\ref{equiv1loimax}), we have that for all $A, \varepsilon >0$ there exist $T$ such that for all $t \geq T$ and $x \geq A$,
\begin{equation}
\label{(B)}
\frac{U(tx)-U(t)}{a(t)} \geq -\frac{1}{\xi} (1+\varepsilon)(1-A^{\xi}).
\end{equation}
Using (\ref{(A)}) and (\ref{(B)}), we conclude the proof for $\xi <
0$. Second, suppose that $\xi=0$ (the proof for $\xi > 0$ is quite similar). Using relation (\ref{equiv1loimax}) we have that for all $A, \varepsilon >0$ there exist $T$ such that for all $t \geq T$ and $x \geq A$,
\[\frac{U(tx)-U(t)}{a(t)} \geq (1+\varepsilon)\ln(A), \]
which concludes the proof. \findemo \\
 \ \\
\noindent {\bf{Proof of Lemma \ref{lem2}}} $-$ Remark that
\[ \frac{U(tx)-U(t)}{U(ty)-U(t)} = 1+z(t,x,y), \]
with
\[ z(t,x,y) = \frac{U(tx)-U(ty)}{U(ty)-U(t)} = {\underbrace{\frac{U(tyx/y)-U(ty)}{a(ty)}}_{z_1(t,x,y)}} \times {\underbrace{\frac{a(ty)}{U(ty)-U(ty/y)}}_{z_2(t,x,y)}}. \]
Since $y \to \infty$ and $x/y \to d$, relation (\ref{equiv1loimax}) implies that $z_1(t,x,y)$ converges to
$\varphi_{\xi}(d)$. Since $y \to \infty$ and $ty \to \infty$, Lemma
\ref{lem1} implies that $z_2(t,x,y)$ converges to $\min(0,\xi)$. Thus,
$1+z(t,x,y) \to 1$ if $\xi \geq 0$ and $1+z(t,x,y) \to c^{\xi}$ if
$\xi < 0$. Remarking that $\varphi_{\xi}(y)/\varphi_{\xi}(x) \to 1$ if $\xi \geq 0$ and $\varphi_{\xi}(y)/\varphi_{\xi}(x)
\to c^{-\xi}$ if $\xi < 0$ concludes the proof. \findemo \\
 \ \\
\noindent {\bf{Proof of Lemma \ref{lemZ}}} $-$ Let $W_1,\ldots,W_n$ be independent standard uniform random variables and $W_{1,n} \leq \ldots \leq W_{n,n}$ the corresponding order statistics.
\begin{itemize}
\item[i)] Remarking that 
\[ (1-K_{i,n}/N_n)_{1 \leq i \leq n} = (F(X_{n-i+1,n}))_{1 \leq i \leq n} \egloi (W_{n-i+1,n})_{1 \leq i \leq n}, \]
it follows that
\begin{eqnarray*}
(U(N_n/K_{i,n}))_{1 \leq i \leq n}=(F^{\leftarrow}(1-K_{i,n}/N_n))_{1 \leq i \leq n}  & \egloi & (F^{\leftarrow}(W_{n-i+1,n}))_{1 \leq i \leq n} \\
 \ & \egloi & (X_{n-i+1,n})_{1 \leq i \leq n}, 
\end{eqnarray*}
which conclude the demonstration.
\item[ii)] We have $K_{k,n} \egloi W_{k,n}/W_{1,n}$. Using R\'{e}nyi
representation \cite{arnbalnag92} p.72, $K_{k,n} \egloi T_k/T_1$ where $T_k$ is the sum of $k$ independent standard exponential random variables. Another use of R\'{e}nyi representation leads to $K_{k,n} \egloi 1/W_{1,k-1}$. Since $W_{1,k-1} \cvps 0$ as $n \to \infty$, we prove that $K_{k,n} \cvps \infty$ as $n \to \infty$. Similarly, $N_n \cvps \infty$ as $n \to \infty$. \\
Now, remark that $K_{k,n}/K_{k',n} \egloi W_{k,n}/W_{k',n}$. R\'{e}nyi representation yields $K_{k,n}/K_{k',n} \egloi T_k/T_{k'} = (k/k') [T_k/k]/[T_{k'}/k'] \cvps c$ as $n \to \infty$. The proof of $K_{k,n}/N_n \cvps 0$ as $n \to \infty$ is similar.\\
Finally, Lemma \ref{lem2} and Lemma \ref{lemZ} i) ii) imply that:
\[ (1+Z_n) \frac{\varphi_{\xi}(1/k')}{\varphi_{\xi}(1/k)} \cvP 1, \]
which concludes the proof using $k/k' \to c$.
\item[iii)] Since $1/K_{k,n} \egloi W_{1,k-1} \egloi 1-W_{k-1,k-1}$, and remarking that
\[ k(1-W_{k-1,k-1}) \cvloi {\rm{Exp}}(1), \]
it follows that $k/K_{k,n}$ has asymptotically a standard exponential distribution.\\
Now, let $E_1, \ldots, E_n$ be independent standard exponential random variables and $E_{1,n} \leq \ldots \leq E_{n,n}$ the corresponding order statistics. Dekkers and de Haan (\cite{dekhaa89}, Lemma 2.1) shows that
\[ \sqrt{\frac{k}{c-1}} (E_{n-k/c+1,n}-E_{n-k+1,n}-\ln(c)) \cvloi {\cal{N}}(0,1).\]
Since $K_{k,n} \egloi \exp(E_{n,n}-E_{n-k+1,n})$ and $K_{k',n} \egloi \exp(E_{n,n}-E_{n-k/c+1,n})$, the $\delta$-method applied to the function $\varphi_{\xi}({\rm{e}}^{-x})$ concludes the proof. \findemo
\end{itemize}
 \ \\
\noindent {\bf{Proof of Lemma \ref{lem3}}} $-$ This proof is inspired
by the one of \cite{res87}, Lemma 0.13. We only give the proof for $\xi >0$ (the case $\xi = 0$ is quite similar). For $j \in \N$, let
\[ \beta_j=\exp \left [ \frac{j}{\xi} \ln(1+\xi) \right ]. \]
Using relation (\ref{equiv1loimax}), we have that for any $\varepsilon
\in ]0,1-(1+\xi)^{-1}[$ there exist $t_0$ such that for all $t \geq
t_0$ and for all $j \in \N \backslash \{0\}$,
\begin{equation}
\label{(C)}
1-\varepsilon \leq \frac{U(t \beta_j)-U(t \beta_{j-1})}{a(t \beta_{j-1})} \leq 1+\varepsilon.
\end{equation}
Furthermore, since $a(.)$ is regularly varying at infinity with index
$\xi$ (see 
\cite{res87}, Proposition 0.8 v) and Proposition 0.12), we have that
for any $\varepsilon \in ]0,1-(1+\xi)^{-1}[$ there exist $t_0$ such
that for all $t \geq t_0$ and for all $j \in \N \backslash \{0\}$,
\begin{equation}
\label{(D)}
\beta_1^{\xi}(1-\varepsilon) = (1+\xi)(1-\varepsilon) \leq \frac{a(t \beta_j)}{a(t \beta_{j-1})} \leq (1+\xi)(1+\varepsilon).
\end{equation}
Let $N \in \N$. We have:
\[ \frac{U(t \beta_{N})-U(t)}{a(t)}=\sum_{j=1}^{N} \frac{U(t \beta_j)-U(t \beta_{j-1})}{a(t \beta_{j-1})} \frac{a(t \beta_{j-1})}{a(t)}. \]
Using (\ref{(C)}), we find that for $t \geq t_0$:
\begin{equation}
\label{(D')}
(1-\varepsilon)\sum_{j=1}^{N}\frac{a(t \beta_{j-1})}{a(t)} \leq \frac{U(t \beta_{N})-U(t)}{a(t)} \leq (1+\varepsilon)\sum_{j=1}^{N}\frac{a(t \beta_{j-1})}{a(t)}.
\end{equation}
Remarking that
\[ \frac{a(t \beta_{j-1})}{a(t)} = \prod_{i=1}^{j-1} \frac{a(t \beta_i)}{a(t \beta_{i-1})}, \]
(\ref{(D)}) and (\ref{(D')}) imply that there exist ${\tilde{\beta_1}}, {\tilde{\beta_2}} > 0$ such that for $t \geq t_0$:
\begin{equation}
\label{(E)}
{\tilde{\beta_1}}[(1+\xi)(1-\varepsilon)]^{N} - {\tilde{\beta_1}} \leq \frac{U(t \beta_{N})-U(t)}{a(t)} \leq {\tilde{\beta_2}}[(1+\xi)(1+\varepsilon)]^{N} - {\tilde{\beta_2}}.
\end{equation}
Let $N_x=[\xi / \ln(1+\xi)]\ln(x)$ (i.e. $x=\exp[N_x/\xi
\ln(1+\xi)]$). Remarking that $\lfloor N_x \rfloor \leq N_x
\leq \lfloor N_x \rfloor +1$ implies that $\beta_{\lfloor N_x
\rfloor } \leq x \leq \beta_{\lfloor N_x \rfloor +1}$. Since $U$ is a non-decreasing function, we find that
\[ \frac{U(t \beta_{\lfloor N_x \rfloor })-U(t)}{a(t)} \leq
\frac{U(tx)-U(t)}{a(t)} \leq \frac{U(t \beta_{\lfloor N_x \rfloor +1})-U(t)}{a(t)}, \]
and, using (\ref{(E)}),
\[ {\tilde{\beta_1}}[(1+\xi)(1-\varepsilon)]^{\lfloor N_x \rfloor} - {\tilde{\beta_1}}
\leq \frac{U(tx)-U(t)}{a(t)} \leq
{\tilde{\beta_2}}[(1+\xi)(1+\varepsilon)]^{\lfloor N_x \rfloor +1} - {\tilde{\beta_2}}. \]
Thus, there exist $\beta_1$ and $\beta_2$ such that:
\begin{equation}
\label{ouf}
\beta_1[(1+\xi)(1-\varepsilon)]^{N_x} - {\tilde{\beta_1}} \leq \frac{U(tx)-U(t)}{a(t)} \leq \beta_2[(1+\xi)(1+\varepsilon)]^{N_x} - {\tilde{\beta_2}}.
\end{equation}
Remarking that
\[ [(1+\xi)(1-\varepsilon)]^{N_x} = x^{\xi-\eta_1}, \]
with $\eta_1 = -\xi \ln(1-\varepsilon)/\ln(1+\xi) \in ]0,\xi[$ and
\[ [(1+\xi)(1+\varepsilon)]^{N_x} = x^{\xi+\eta_2}, \]
with $\eta_2 = \xi \ln(1+\varepsilon)/\ln(1+\xi) > 0$ concludes the proof. \findemo \\
 \ \\
\noindent {\bf{Proof of Lemma \ref{lemcvloiZ}}} $-$ The first step of the proof consists in establishing the following expansion:
\begin{equation}
\label{resulinter}
-\frac{\varphi_{\xi}(1/c)}{Z_n} = \varphi_{\xi}(K_{k',n}) \left [ 1 - \frac{\sigma}{\sqrt{k} \varphi_{\xi}(K_{k',n}/K_{k,n})}Y_n \right ] \left [ 1 + o_{{\rm{P}}} \left ( \frac{1}{\varphi_{\delta}(k)} \right ) \right ].
\end{equation}
Let $V(t)=U({\rm{e}}^t)$. From Lemma \ref{lemZ} i), we have:
\[ Z_n \egloi \underbrace{\frac{V'[\ln(N_n/K_{k',n})]}{U(N_n/K_{k',n})-U(N_n)}}_{Z_{1,n}} \times \underbrace{\frac{U(N_n/K_{k,n})-U(N_n/K_{k',n})}{V'[\ln(N_n/K_{k',n})]}}_{Z_{2,n}}.\]
Clearly,
\[ -\frac{1}{Z_{1,n}} = \int_0^{\ln(K_{k',n})} \frac{V'[\ln(N_n/K_{k',n})+s]}{V'[\ln(N_n/K_{k',n})]} ds, \]
and conditions $(H1)$ and $(H2)$ imply that uniformly on $[0,\ln(K_{k',n})]$,
\[ \frac{V'[\ln(N_n/K_{k',n})+s]}{V'[\ln(N_n/K_{k',n})]} = {\rm{e}}^{\xi s} \left [ 1 + o_{{\rm{P}}} \left ( \frac{1}{\varphi_{\delta}(k)} \right ) \right ]. \]
Thus,
\begin{equation}
\label{Z_{1,n}p}
-\frac{1}{Z_{1,n}} =\left [ 1 + o_{{\rm{P}}} \left ( \frac{1}{\varphi_{\delta}(k)} \right ) \right ] \int_0^{\ln(K_{k',n})} {\rm{e}}^{\xi s} ds = \varphi_{\xi}(K_{k',n}) \left [ 1 + o_{{\rm{P}}} \left ( \frac{1}{\varphi_{\delta}(k)} \right ) \right ].
\end{equation}
The proof of
\begin{equation}
\label{Z_{2,n}p}
Z_{2,n} =\varphi_{\xi}(K_{k',n}/K_{k,n}) \left [ 1 + o_{{\rm{P}}} \left ( \frac{1}{\varphi_{\delta}(k)} \right ) \right ],
\end{equation}
follows the same lines. Collecting (\ref{Z_{1,n}p}) and (\ref{Z_{2,n}p}) yields
\[ -\frac{\varphi_{\xi}(K_{k',n}/K_{k,n})}{Z_n} = \varphi_{\xi}(K_{k',n}) \left [ 1 + o_{{\rm{P}}} \left ( \frac{1}{\varphi_{\delta}(k)} \right ) \right ], \]
which proves (\ref{resulinter}) by remarking that
\[ \frac{\varphi_{\xi}(1/c)}{\varphi_{\xi}(K_{k',n}/K_{k,n})} = 1 - \frac{\sigma}{\sqrt{k} \varphi_{\xi}(K_{k',n}/K_{k,n})} Y_n. \]
Now, remark that (\ref{resulinter}) can be rewritten $v_n=u_n(1+\varepsilon_n)$ with $u_n=\varphi_{\xi}(K_{k',n})$, $v_n=-\varphi_{\xi}(1/c)/Z_n$ and
\[ \varepsilon_n = \frac{-\sigma}{\sqrt{k} \varphi_{\xi}(K_{k',n}/K_{k,n})} Y_n+o_{{\rm{P}}} \left ( \frac{1}{\varphi_{\delta}(k)} \right ) = O_{{\rm{P}}} \left ( k^{-1/2} \right ) + o_{{\rm{P}}} \left ( \frac{1}{\varphi_{\delta}(k)} \right ) \cvP 0, \]
from Lemma \ref{lemZ} iii). The second step of the proof is dedicated to the study of $\varphi_{\xi}^{*} \left ( - \varphi_{\xi}(1/c)/Z_n \right )$. Five cases have to be considered:\\
\noindent {\underline{If $\xi >0$}}, since $u_n \cvP +\infty$ as $n \to \infty$ (Lemma \ref{lemZ} iii)), Lemma \ref{lemphi} III) entails that $\varphi_{\xi}^{*}(u_n) \eqproba \varphi_{\xi}^{*}(v_n)$, i.e.,
\begin{equation}
\label{cas1}
\varphi_{\xi}^{*} \left ( - \frac{\varphi_{\xi}(1/c)}{Z_n} \right ) \eqproba \varphi_{\xi}^{*}[\varphi_{\xi}(K_{k',n})] = (K_{k',n})^{\xi+\ind\{\xi=0\}},
\end{equation}
by Lemma \ref{lemphi} I).\\
{\underline{If $\xi = 0$}}, we have $u_n \cvP +\infty$ as $n \to \infty$ and $u_n-v_n=u_n \varepsilon_n=o_{{\rm{P}}}[\ln(K_{k',n})\ln(k)]=o_{{\rm{P}}}(1)$. Thus, Lemma \ref{lemphi} II) implies (\ref{cas1}).\\
{\underline{If $-1/2<\xi<0$}}, we have from Lemma \ref{lemZ} ii) that
$u_n \cvP -1/\xi$ as $n \to \infty$. Remarking that
$\varepsilon_n=o_{{\rm{P}}}[1/\varphi_{\delta}(k)] =
o_{{\rm{P}}}(k^{\xi})$ implies
\[ \frac{\varepsilon_n}{\varphi_{\xi}^{*}(u_n)} =o_{{\rm{P}}} [(k/K_{k',n})^{\xi}] = o_{{\rm{P}}}(1), \]
which entails (\ref{cas1}) by Lemma \ref{lemphi} IV) i).\\
{\underline{If $\xi=-1/2$}}, we have $u_n \cvP -1/\xi$ as $n \to
\infty$. Remarking that $\varepsilon_n \eqproba -\sigma/[\sqrt{k}
\varphi_{\xi}(1/c)] Y_n$ yields:
\[ \frac{\varepsilon_n}{\varphi_{\xi}^{*}(u_n)} \eqproba -\frac{\sigma c^{-1/2}}{\varphi_{\xi}(1/c)} Y_n \sqrt{ \frac{K_{k',n}}{k'} } = \alpha_n, \]
where $\alpha_n$ does not converges in probability to $\infty$ or $1$ as $n \to \infty$ (see Lemma \ref{lemZ} iii)). Thus, from Lemma \ref{lemphi} IV) i), we have $\varphi_{\xi}^{*}(u_n) \eqproba \varphi_{\xi}^{*}(v_n)(1-\alpha_n)$ i.e.
\begin{equation}
\label{cas2}
\varphi_{\xi}^{*} \left ( - \frac{\varphi_{\xi}(1/c)}{Z_n} \right ) \eqproba (K_{k',n})^{-1/2} \left ( 1 + \frac{\sigma c^{-1/2}}{\varphi_{\xi}(1/c)} Y_n\sqrt{ \frac{K_{k',n}}{k'} } \right ) = k^{-1/2}T_n.
\end{equation}
{\underline{If $\xi<-1/2$}}, we have $u_n \cvP -1/\xi$ as $n \to \infty$, $\varepsilon_n \eqproba -\sigma/[\sqrt{k} \varphi_{\xi}(1/c)] Y_n$ and $\varepsilon_n/\varphi_{\xi}^{*}(u_n) \cvP \infty$ as $n \to \infty$. Thus, Lemma \ref{lemphi} IV) ii) implies that $\varphi_{\xi}^{*}(v_n) \eqproba -\varepsilon_n$ i.e.
\begin{equation}
\label{cas3}
\varphi_{\xi}^{*} \left ( - \frac{\varphi_{\xi}(1/c)}{Z_n} \right ) \eqproba \frac{\sigma k^{-1/2}}{\varphi_{\xi}(1/c)} Y_n.
\end{equation}
Lemma \ref{lemZ} iii) and (\ref{cas1})-(\ref{cas3}) conclude the proof. \findemo

\bibliography{mabiblio}

\end{document}